\documentclass[12pt, reqno]{amsart}
\usepackage{amsmath,amssymb,amsthm}
\usepackage{amsmath, amssymb}
\usepackage{amsthm, amsfonts, mathrsfs}
\usepackage{mathptmx}
\usepackage{fullpage}
\usepackage{amsfonts,graphicx}
\usepackage{dutchcal}
\usepackage{bm}
\numberwithin{equation}{section}
\usepackage[colorlinks=true, pdfstartview=FitV, linkcolor=blue, citecolor=blue, urlcolor=blue]{hyperref}

\newtheorem {theorem}{Theorem}[section]
\newtheorem {lemma}[theorem]{{\bf Lemma}}

\newtheorem {proposition}[theorem]{{\bf Proposition}}
\theoremstyle{remark}
\newtheorem {remark}{{\bf Remark}}[section]

\theoremstyle{plain} \numberwithin {equation}{section}

\newcommand{\R}{{\mathbb R}}
\def\div{ \hbox{\rm div}\,  }

\def\nn{\nonumber}

\def\f{\frac}

\def\u{ \mathbf{u} }
\def\v{ \mathbf{v} }
\def\T{ \mathbb{T} }

\begin{document}
\title{  Stability Threshold of the 2D Boussinesq System Near Couette\\[1ex] Flow in an Infinite  Channel  }

\author{Tao Liang}

\address[T. Liang]{ School of Mathematics,
	South China University of Technology,
	Guangzhou, 510640, China}
\email{taolmath@163.com}

\author[J. Wu]{Jiahong Wu}
\address[J. Wu]{Department of Mathematics, University of Notre
	Dame, Notre Dame, IN 46556, USA } \email{jwu29@nd.edu}

 \author[X. Zhai]{Xiaoping Zhai}
\address[X. Zhai]{School of Mathematics and Statistics, Guangdong University of Technology,
Guangzhou, 510520, China} \email{pingxiaozhai@163.com (Corresponding author)}

\footnote{\today}



\begin{abstract}
In this paper, we study the stability threshold of the two-dimensional Boussinesq equations around the Couette flow in an infinite channel $\mathbb{R} \times [-1, 1]$ under no-slip boundary conditions. We prove that the Couette flow is asymptotically stable under initial perturbations satisfying
$\| \v^{\mathrm{in}} -(y,0)\|_{H^2} \le \varepsilon_0 \nu^{\frac12}$, and $\| \rho^{\mathrm{in}}-1 \|_{H^1} + \big\| |\partial_x|^{\frac13} \rho^{\mathrm{in}} \big\|_{H^1} \le \varepsilon_1 \nu^{\frac56}.
$
Compared with the work of Masmoudi, Zhai, and Zhao [J. Funct. Anal., 284 (2023), 109736], where the asymptotic stability of the 2D Navier-Stokes-Boussinesq system around Couette flow in a finite channel $\mathbb{T} \times [-1, 1]$ was established, our result improves the stability threshold for the temperature from $\nu^{\frac{11}{12}}$ to $\nu^{\frac56}$.
\end{abstract}

\maketitle

\tableofcontents
\section{Introduction and the main result}
\subsection {Two dimensional Navier-Stokes Boussinesq equations}

In this paper, we study the two dimensional Navier-Stokes-Boussinesq system under the no slip boundary condition in an infinite channel domain $(x,y) \in \R \times [-1, 1]$, which is formulated as follows:
\begin{eqnarray}\label{quanhs}
\left\{\begin{aligned}
&\partial_t \v+ \v \cdot\nabla \v - \nu (  \partial_x^2 + \partial_y^2) \v  + \nabla P  = \rho g \mathbf{e}_2  ,\\
& \partial_t \rho + \v \cdot\nabla\rho - \mu (  \partial_x^2 + \partial_y^2) \rho  = 0,\\
& \div \v = \partial_x v_1 + \partial_y v_2 = 0,\\
& \v(t, x, \pm 1) = (\pm 1, 0),   \quad \rho(t, x, \pm 1) = c_0,\\
& \v(0, x,y) = \v^{\mathrm{in}}(x,y),\qquad \rho(0,x,y) = \rho^{\mathrm{in}}(x, y).
\end{aligned}\right.
\end{eqnarray}
Here,
$\v=(v^1,v^2)$, $\rho$, and $P$ denote the velocity, temperature, and pressure, respectively. The constants
$\nu$, $\mu$, and $g$ represent the viscosity coefficient, thermal diffusivity coefficient, and gravitational constant, respectively, while $\mathbf{e}_2 = (0,1)$ is the unit vector in the vertical direction. For simplicity and without loss of generality, we normalize the constants by setting
   $g = c_0 = 1$.

This paper focuses on the stability problem near the following steady-state solution:
\begin{align*}
	\v_s = ( y, 0), \quad  \rho_s = 1, \quad P_s = y + c.
\end{align*}

Let $ \u = \v - \v_s$, $\theta = \rho - \rho_s$, and $p = P - P_s$. Then the system \eqref{quanhs} can be reformulated as the following equations for $(\u, \theta, p)$:
\begin{eqnarray}\label{rewrite}
\left\{\begin{aligned}
&\partial_t \u+ y\partial_x \u     + \u \cdot\nabla \u - \nu (  \partial_x^2 + \partial_y^2) \u  + \nabla p  =  \binom{0}{ \theta},\\
& \partial_t \theta + y\partial_x \theta + \u \cdot\nabla\theta - \mu (  \partial_x^2 + \partial_y^2) \theta  = 0,\\
& \div \u =  0,\\
& \u(t, x, \pm 1) = 0,   \quad \theta(t, x, \pm 1) = 0.
\end{aligned}\right.
\end{eqnarray}
By introducing the vorticity $ \omega = \nabla \times \u = \partial_y u^1 - \partial_x u^2$, the above system \eqref{rewrite} can be reformulated as follows
\begin{eqnarray}\label{rewrite1}
\left\{\begin{aligned}
&\partial_t \omega+ y\partial_x \omega   + \u \cdot\nabla \omega - \nu (  \partial_x^2 + \partial_y^2) \omega    =  -\partial_x \theta,\\
& \partial_t \theta + y\partial_x \theta + \u \cdot\nabla\theta - \mu (  \partial_x^2 + \partial_y^2) \theta  = 0,\\
& \u = \nabla^{\perp} \psi = (\partial_y \psi,-\partial_x \psi ), \quad \Delta \psi =\omega,\\
& \psi(t, x, \pm 1) = \partial_{y} \psi(t, x, \pm 1) = 0, \quad  \theta(t, x, \pm 1) = 0,\\
& \omega(0, x, y) = \omega^{\mathrm{in}}(x,y), \quad \theta(0, x, y) = \theta^{\mathrm{in}}(x,y).
\end{aligned}\right.
\end{eqnarray}
 The transition from laminar to turbulent flow has been a fundamental challenge in fluid dynamics since Reynolds' pioneering experiments \cite{Re}. Although certain laminar flows remain linearly stable at all Reynolds numbers \cite{DR,Rom}, they can exhibit nonlinear instability when subjected to finite amplitude perturbations at high Reynolds numbers \cite{DHB,TA}. Kelvin first observed \cite{Kel} that the basin of attraction of laminar flow shrinks as $Re \to \infty$, making nonlinear instability possible. This observation motivates the central question posed by Trefethen et al. \cite{T}: determining the transition threshold the minimum disturbance amplitude that triggers instability and its scaling behavior with the Reynolds number.

To address this question, we adopt the stability threshold framework introduced by Bedrossian, Germain, and Masmoudi \cite{BGM-BAMS}. In this framework, one seeks exponents $\alpha=\alpha(Y_1,Y_2)$ and $\beta=\beta(Y_1,Y_2)$ associated with given norms  $\|\cdot\|_{Y_1}$ and $\|\cdot\|_{Y_2}$ such that:
\begin{itemize}
  \item  $\left\|\omega^{\mathrm{in}}\right\|_{Y_1} \leq \nu^\alpha$ \hbox{ and }$\left\|\theta^{\mathrm{in}}\right\|_{Y_2} \leq \mu^\beta \Rightarrow$  \emph{stable};\\
  \item $\left\|\omega^{\mathrm{in}}\right\|_{Y_1} \gg \nu^\alpha$ \text { or }$\left\|\theta^{\mathrm{in}}\right\|_{Y_2} \gg \mu^\beta \Rightarrow$  \emph{unstable}.
\end{itemize}
The exponents $\alpha $ and $\beta$ are referred to as the transition threshold.

\subsection  {Historical comments}
The Navier-Stokes-Boussinesq (NSB) system serves as a fundamental model for buoyancy driven flows, with broad applications across atmospheric science, oceanic circulation, mantle convection, and engineering heat transfer \cite{Gill,Majda}. Based on the Boussinesq approximation, this approach retains density variations exclusively in the buoyancy term such as those caused by temperature or salinity differences while neglecting them in inertial terms. This model effectively captures essential phenomena including Rayleigh-B\'enard convection and stratified shear flows \cite{Constantin,Getling}.

Before presenting our main results, we provide a brief review of existing studies on the stability of perturbations about steady states for System \eqref{quanhs}. The stability analysis is conventionally classified into three categories according to the domain type: perturbations defined on $\mathbb{T} \times \mathbb{R}$, on $\mathbb{T} \times [-1,1]$, and on $\mathbb{R} \times \mathbb{R}$.
\vskip .1in
$(\mathrm{I})$ The domain $\T \times \R$.
\begin{itemize}
  \item
  For perturbations about the shear flow $\mathbf{v}_s = (U(y), 0)$ in the absence of hydrostatic equilibrium ($\rho_s = 1$), Tao and Wu \cite{Tao2017JDE} established the linear stability of System \eqref{rewrite1} in the upper half-space $\mathbb{T} \times \mathbb{R}^+$, demonstrating that enhanced dissipation occurs even when dissipation is only present in the vertical direction.
 Subsequently, Deng et al. \cite{Deng2021JFA} proved the asymptotic stability of the steady state under the condition that the initial perturbations $(\omega^{\mathrm{in}}, \theta^{\mathrm{in}})$ satisfy
    \begin{equation*} 
        \left\|\omega^{\mathrm{in}}\right\|_{H^s} \leq c \nu^{\alpha}, \quad    \left\|\theta^{\mathrm{in}}\right\|_{H^s} \leq c \nu^{\beta}, \quad \left\|\left|D_x\right|^{\frac{1}{3}} \theta^{\mathrm{in}}\right\|_{H^s} \leq c \nu^{\delta},
    \end{equation*}
    where $s>1$, $\alpha\ge\frac12$, $\delta\ge \alpha+\frac13$, $\beta\ge\delta-\alpha+\frac23.$
    The threshold of \cite{Deng2021JFA}  was improved by Zhang and Zi \cite{Zhang2023JMPA} to
	$$
	\left\|\u^{\mathrm{in}}\right\|_{H^{s+1}} \leq \varepsilon_0 \nu^{\frac{1}{3}},\quad \left\|\theta^{\mathrm{in}}\right\|_{H^s}+\nu^{\frac{1}{6}}\left\|\left|D_x\right|^{\frac{1}{3}} \theta^{\mathrm{in}}\right\|_{H^s} \leq \varepsilon_1 \nu^{\frac{5}{6}},
	$$
	where $s>7$. Later, Niu and Zhao \cite{Niu2024ARxiv}  combined the quasi-linearization technique introduced in \cite{Chenqi2020ARMA,Chenqi2024AMS} with the time-dependent elliptic operator
$ \Lambda_t^2 = 1-\partial_x^2 - (\partial_y + t\partial_x)^2$,
constructed in \cite{Deng2021JFA,Wei2023Tunisi}, and applied Fourier analysis to further improve the stability results from \cite{Zhang2023JMPA}. Specifically, they showed that the steady state is asymptotically stable provided the initial perturbations satisfy
$$
	\left\|\u^{\mathrm{in}}\right\|_{H^{s+1}} \leq \epsilon_0 \nu^{\frac{1}{3}}, \quad \big\|\left\langle\partial_x\right\rangle \theta^{\mathrm{in}}\big\|_{H^s} \leq \epsilon_0 \nu^{\frac{2}{3}},
	$$
	where $s>5$.

 \vskip .1in
  \item
  When examining perturbations around the shear flow $\mathbf{v}_s = (U(y), 0)$ in the presence of hydrostatic equilibrium ($\rho_s = \rho_s(y)$), Yang and Lin \cite{Lin2018JMFM} established the linear stability of system \eqref{quanhs} for the inviscid Couette flow $\mathbf{v}_s = (y, 0)$ with an exponentially stratified density $\rho_s(y)$. This result was later extended by Bianchini et al. \cite{Dolce2021India} to near-Couette flows $\mathbf{v_s} = (U(y), 0)$ satisfying $ \| U(y) - y\|_{H^s} \ll 1$. Subsequently, Bedrossian et al. \cite{Dolce2023CPAM} demonstrated shear-buoyancy instability for Gevrey-class initial perturbations of size $\varepsilon$, under specific constraints on the Richardson number. In the non-diffusive case, Masmoudi et al. \cite{Zhao2022ARMA} established stability for Gevrey$-\frac{1}{s}$ initial data with $\frac{1}{3} < s \leq 1$. For the fully dissipative system, Zillinger obtained nonlinear stability and enhanced dissipation in \cite{Zillinger1}, and later extended these results in \cite{Zillinger2} to oscillatory temperature profiles under vertical-only viscous and thermal dissipation.
  Additionally, Zhai and Zhao \cite{Zhao2023SIAM} derived a Sobolev stability threshold of $\nu^{\frac{1}{2}}$ for the two-dimensional Boussinesq equations with $\mu = \nu$ around Couette flow, under the condition that the Richardson number satisfies $\gamma^2 > \frac{1}{4}$. This threshold was later improved by Knobel \cite{Knobel} to $\nu^{\frac{1}{3}}$ using a method that reduces nonlinear stability to a linear analysis. Concurrently, Ren and Wei \cite{renxiaoxia2026} weakened the regularity assumptions required in \cite{Zhao2023SIAM}. For the case of distinct viscosity $\nu$ and thermal diffusivity $\mu$ with $\gamma^2 > \frac{1}{4}$,  Ren and Wei \cite{renxiaoxia2026} established asymptotic stability under the condition $\frac{\nu + \mu}{2\gamma\sqrt{\nu \mu}} < 2 - \varepsilon$ with $0 < \varepsilon < \frac{1}{\gamma}$, assuming the initial data satisfy
  $$
\big\|\v^{\mathrm{in}}-(y,0)\big\|_{H^{s+\frac12}}+\big\|\rho^{\mathrm{in}}+\gamma^2 y-1\big\|_{H^{s+\frac12}}\leq c \min\{\nu, \mu
\}^{\frac{1}{3}},
$$
where  $s >\frac32$. In a separate work,
Liang et al. \cite{liangtaoscm} proved asymptotic stability in the regime $\mu^3 \leq \nu \leq \mu^{\frac{1}{3}}$, provided the initial perturbations $(\mathbf{u}^{\mathrm{in}}, \theta^{\mathrm{in}})$ satisfy
  $$
\qquad\big\|\v^{\mathrm{in}}-(y,0)\big\|_{H^{N+1}}+\big\|\rho^{\mathrm{in}}+\gamma^2 y-1\big\|_{H^{N+2}}\leq c \min\{\nu, \mu
\}^{\frac{1}{2} + \frac{4}{3}\alpha},$$
where $N \ge 6$ and $ 0 < \alpha \ll 1$.
  In subsequent work, Zelati, Zotto, and Widmayer \cite{ZelatiARxiv} investigated the stability threshold for system \eqref{quanhs} in the three-dimensional setting. Their analysis suggests, from a mathematical perspective, that while a Richardson number $\gamma^2 > \frac{1}{4}$ is sufficient for stability, it may not be a necessary condition in the 3D case.
\end{itemize}

\vskip .1in
$(\mathrm{II})$ The domain $\T \times [-1,1]$.
\begin{itemize}
	\item
When studying perturbations around the shear flow $\mathbf{v}_s = (U(y), 0)$ in the presence of hydrostatic equilibrium ($\rho_s = 1$), Masmoudi et al. \cite{Masmoudi2023JFA} established the nonlinear stability of system \eqref{rewrite1} for the Couette flow $U(y) = y$  if the initial velocity and initial temperature $(\v^{\mathrm{in}},\rho^{\mathrm{in}})$ satisfies $$\|\v^{\mathrm{in}}-(y,0)\|_{H_{x,y}^2}\leq \varepsilon_0 \min\{\nu,\mu\}^{\f12},\quad \|\rho^{\mathrm{in}}-1\|_{H_x^{1}L_y^2}\leq \varepsilon_1 \min\{\nu,\mu\}^{\f{11}{12}}.$$
Liang et al. \cite{liangtaojimj} investigated stability around near-Couette shear flows $(U(y), 0)$ under Navier slip boundary conditions. They proved nonlinear stability under the condition that the initial vorticity lies in an anisotropic Sobolev space with norm of order $O(\min { \mu^{\frac{1}{2}}, \nu^{\frac{1}{2}}})$, while the initial temperature perturbation lies in an anisotropic Sobolev space of size $O(\min { \mu, \nu})$.
\end{itemize}

\vskip .1in
$(\mathrm{III})$ The domain $\R \times\R$.
\begin{itemize}
	\item In the whole plane,
Wang and Wang \cite{wangwang} established the asymptotic stability of the steady state under the condition that the initial perturbations $(\u^{\mathrm{in}}, \theta^{\mathrm{in}})$ satisfy
	$$
	\left\|\omega^{\mathrm{in}}\right\|_{H^1 \cap L^1} \leq \epsilon_0 \nu^{\frac{3}{4}}, \quad\left\|\theta^{\mathrm{in}}\right\|_{H^1 \cap L^1} \leq \epsilon_1 \nu^{\frac{5}{4}}.
	$$
In a related work, Chen et al. \cite{chenwangyang} derived the asymptotic stability threshold for the 2D Boussinesq system with Couette flow under the assumption that the initial perturbations lie in Sobolev spaces with controlled low horizontal frequencies. They showed that the stability threshold is at most $\left\{\frac{1}{3}+, \frac{2}{3}+\right\}$.
Further extending these results, Liang et al. \cite{liangtaofenshujie} proved the nonlinear asymptotic stability of the Couette flow for the two-dimensional Boussinesq equations with fractional dissipation $(-\Delta)^{\alpha}$, where $\alpha \in \left(\tfrac{\sqrt{2}}{2}, 1\right]$. They obtained a stability threshold of ${\nu}^{\frac{1}{2\alpha}}$.
\end{itemize}

\vskip .1in
$(\mathrm{IV})$
For the 3D Boussinesq equations via the Couette flow, we refer to \cite{ZZ2025,ZZW2024,CWW202501, ZelatiARxiv} and the reference therein.

Moreover, we refer to \cite{Cao2012, Chae2006,  Danchin2011,Houyizhao, Lai2011, Lijinkai2016} for results concerning the well-posedness of system \eqref{quanhs} in both the whole space and bounded domains.

\subsection{Main result}
To make the statement precise, we define, for $(x,y) \in \R \times [-1, 1]$,
$$ f_k(t,y) =   \int_{\R} f(t,x,y)e^{-ikx} \mathrm{d}x.$$
By taking the Fourier transform with respect to the horizontal variable $x$, we decompose $(\omega, \theta)$ into their Fourier modes corresponding to each wave number $k \in \R$ as follows:
\begin{eqnarray}\label{equation_k}
\left\{\begin{aligned}
&\partial_t \omega_k+ ik y \omega_k   - \nu (  \partial_y^2 -  k^2  ) \omega_k    =  -i k \theta_k - (\u \cdot \nabla \omega)_k ,\\
& \partial_t \theta_k + i k y \theta_k  - \nu (  \partial_y^2 -  k^2) \theta_k  = - (\u \cdot \nabla \theta)_k,\\
&  \u_k  = (\partial_y \psi_k,-ik \psi_k ) , \quad  \Delta_{k} \psi_k =(-k^2 + \partial_y^2)\psi_k=\omega_k,\\
& \psi_k (t, \pm 1) = \partial_{y} \psi_k(t, \pm 1) = 0, \quad  \theta_k(t, \pm 1) = 0.
\end{aligned}\right.
\end{eqnarray}

\begin{theorem}\label{thm1.2}
	Consider the case $\nu = \mu$. Let $(\u, \theta)$ be a solution to the system \eqref{rewrite} with initial data $(\u^{\mathrm{in}}, \theta^{\mathrm{in}})$. Then there exist positive constants $\nu_0$, $\varepsilon_0$, $\varepsilon_1$, and $C > 0$, independent of $\nu$, such that if $0 < \nu \le \nu_0$ and the initial data satisfy the smallness condition
	\begin{align*}
		& \| \u^{\mathrm{in}} \|_{H^2} \le \varepsilon_0 \nu^{\frac{1}{2} },\\
&\| \theta^{\mathrm{in}} \|_{H^1} + \big\| \left| \partial_x\right|^{\frac{1}{3}} \theta^{\mathrm{in}} \big\|_{H^1} \le \varepsilon_1 \nu^{\frac{5}{6} },
	\end{align*}
	then the system \eqref{rewrite} admits a global-in-time solution $(\u, \theta)$. Moreover, the solution satisfies the following uniform stability estimates:
	\begin{align*}
		\big\| E_k[\omega_k]\big\|_{L^1_k(\R)} \le C \varepsilon_0 \nu^{\frac{1}{2}}  \quad \text{and} \quad  \big\| E_k[\theta_k]\big\|_{L^1_k(\R)} \le C \varepsilon_1 \nu^{\frac{5}{6}},
	\end{align*}
	where
	\begin{align*}
		E_k[\omega_k] = \begin{cases} \| \omega_k\|_{L^\infty_t L^2_y} + \nu^{\frac{1}{2} }  \| \omega_k\|_{L^2_t L^2_y} +  \| \u_k\|_{L^\infty_t L^\infty_y} + \nu^{\frac{1}{2} }  \| \u_k\|_{L^2_t L^2_y} , & 0 \le |k| \le 10 \nu ; \\
 \|(1-|y|)^{\frac{1}{2}} \omega_k\|_{L^\infty_t L^2_y} + \nu^{\frac{1}{4} } |k|^{-\frac{1}{4} }  \| \omega_k\|_{L^\infty_t L^2_y} + \nu^{\frac{1}{4} } |k|^{\frac{1}{4} } \| \omega_k\|_{L^2_t L^2_y} \\
  \qquad\qquad\qquad\qquad\quad+   \| \u_k\|_{L^\infty_t L^\infty_y}+ |k|^{\frac{1}{2} }  \| \u_k\|_{L^2_t L^2_y}, & 10 \nu \le |k| \leq 1;\\
   \|(1-|y|)^{\frac{1}{2}} \omega_k\|_{L^\infty_t L^2_y} +   \nu^{\frac{1}{4} } |k|^{\frac{1}{2} } \| \omega_k\|_{L^2_t L^2_y} +  |k|^{\frac{1}{2}} \| \u_k\|_{L^\infty_t L^\infty_y}  + |k|  \| \u_k\|_{L^2_t L^2_y}, &   |k| \geq 1,\end{cases}
	\end{align*}
	and
	\begin{align*}
	E_k[\theta_k] = \begin{cases} \| \theta_k\|_{L^\infty_t L^2_y} + \nu^{\frac{1}{2} }  \| \theta_k\|_{L^2_t L^2_y}  , & 0 \le |k| \le 10 \nu ; \\    \| \theta_k\|_{L^\infty_t L^2_y} + \nu^{\frac{1}{6} } |k|^{\frac{1}{3} } \| \theta_k\|_{L^2_t L^2_y} , & 10 \nu \le |k| \leq 1;  \\ |k|^{\frac{1}{3}} \| \theta_k\|_{L^\infty_t L^2_y} +   \nu^{\frac{1}{6} } |k|^{\frac{2}{3} } \| \theta_k\|_{L^2_t L^2_y} , &   |k| \geq 1.\end{cases}
	\end{align*}
\end{theorem}

\begin{remark}
To the best of our knowledge, this study establishes the first rigorous asymptotic stability result for the two-dimensional Boussinesq system in the vicinity of Couette flow within an infinite channel domain.
\end{remark}

\begin{remark}
	In the work of Masmoudi, Zhai, and Zhao \cite{Masmoudi2023JFA}, the thresholds for $(\omega, \theta)$ on $\mathbb{T} \times [-1,1]$ were obtained as $\left(O(\nu^{\frac12}), O(\nu^{\frac{11}{12}})\right)$. Owing to the one-derivative loss in the buoyancy term, they ingeniously compensated for most of this loss by exploiting the enhanced dissipation of the vorticity, $\nu^{\frac{1}{4}} |k|^{\frac{1}{4}} \|\omega_k\|_{L^2_t L^2_y}$, and that of the temperature, $\nu^{\frac{1}{6}} |k|^{\frac{1}{3}} \|\theta_k\|_{L^2_t L^2_y}$. The fractional operator $|\partial_x|^{\frac16}$ was then used to measure and control the remaining unabsorbed portion of the derivative loss.
	In contrast, in the present paper, it suffices to assume an additional $\frac13$ derivative of regularity on the high-frequency component of $\theta$, which allows us to larger the temperature threshold to $O(\nu^{\frac56})$.
\end{remark}

\begin{remark}
	 The result of Theorem \ref{thm1.2} remains valid in a periodic channel $\mathbb{T} \times [-1,1]$, where, for any fixed wavenumber $k$, the shear-enhanced dissipation rate scales as $\nu^{\frac13}|k|^{\frac23}$. However, in the infinite channel $\mathbb{R} \times [-1,1]$ with continuous wavenumber spectrum, the long-wave limit $k \to 0$ presents a fundamental obstruction: the enhanced dissipation rate $\nu^{\frac13} |k|^{\frac23}$ vanishes, leading to an extremely weak enhanced dissipation effect. Consequently, the dissipation time scale $\sim \nu^{-1/3} |k|^{-\frac23}$ becomes very long. Similarly, the decay rate of inviscid damping, which scales as $\sim \frac{1}{|k| t^2}$, also becomes slow. This permits persistent nonlinear excitation of nearly dissipation-free long-wave modes, potentially preventing decay to equilibrium and inducing instability.
\end{remark}

\begin{remark}
	To control the buoyancy term $\partial_x \theta$, one can normally use the enhanced dissipation effects of $\omega$ and $\theta$ together with the viscous effects. This typically leads to the stability threshold for $\theta$ being lower than that of $\omega$ by $\nu^{\frac12}$. More precisely, we have
	\begin{align*}
		\operatorname{Re} \left\langle \omega_k,  \partial_x \theta \right\rangle &\le \varepsilon\nu^{\frac13}|k|^{\frac23} \| \omega_k\|^2_{L^2_y} + \frac{1}{4\varepsilon} \nu^{-\frac13}|k|^{\frac43} \| \theta_k\|^2_{L^2_y} \\
		&\le \varepsilon\nu^{\frac13}|k|^{\frac23} \| \omega_k\|^2_{L^2_y} + \frac{1}{4\varepsilon} \nu^{-1}( \nu^{\frac13}|k|^{\frac23} \| \theta_k\|^2_{L^2_y} + \nu |k|^2 \| \theta_k\|^2_{L^2_y} ).
	\end{align*}
	If there were no derivative loss in the buoyancy term, the stability threshold for the temperature would be lower than that of the vorticity by $\nu^{\frac13}$. Therefore, we expect that when there is a derivative loss in the buoyancy term, the temperature's stability threshold behaves similarly. To raise the temperature threshold $\nu^{\alpha}$, we need to assume an extra $\frac13$ derivative of regularity on the initial temperature. This compensates for the one-derivative loss from the buoyancy term, allowing the enhanced dissipation and viscous effects to jointly control it. As a result, the threshold for the temperature remains lower than that of the vorticity by $\nu^{\frac13}$, which is the same as in the case without derivative loss.
\end{remark}

\vskip .1in
\subsection*{Notations}
Throughout this paper, we use $ C > 0 $ to denote a generic constant that is independent of all relevant quantities and may vary from line to line. For brevity, we write $ f \lesssim g $ if there exists such a constant $C $ satisfying $  f \leq C g$. In addition, we denote by $ \langle f, g \rangle$ the $ L^2_y([-1, 1])$ inner product of $ f$  and $g $.

\section{Space-time estimates for the linearized Boussinesq equations}\label{sec2}
In this section, inspired by the works \cite{Chenqi2020ARMA, Miao2025Arxiv} on the two-dimensional incompressible Navier--Stokes equations, they establish refined resolvent estimates for both the no-slip and Navier-slip boundary conditions. In particular, by exploiting properties of the Airy function, they achieve precise control over the boundary layers. As a first step, space-time estimates for the vorticity are obtained. For the temperature equation, by decomposing the problem into homogeneous and nonhomogeneous components and treating the enhanced dissipation, boundary layer effects, and the influence of the initial data separately, we derive space-time estimates for the temperature at different horizontal frequencies.
For the nonlinear terms in the first two equations of system \eqref{equation_k}, we exploit the incompressibility condition $\div \u = 0$ to represent their $k$-th Fourier modes as
\begin{align*}
	&- (\u \cdot \nabla \omega)_k = -ik \int_{\R} u^1_\ell \cdot \omega_{k-\ell}  \mathrm{d} \ell - \partial_{y} \int_{\R} u^2_\ell \cdot \omega_{k-\ell}  \mathrm{d} \ell := -ik f^1_{k} - \partial_{y} f^2_{k},\\
	&- (\u \cdot \nabla \theta)_k = -ik \int_{\R} u^1_\ell \cdot \theta_{k-\ell}  \mathrm{d} \ell - \partial_{y} \int_{\R} u^2_\ell \cdot \theta_{k-\ell}  \mathrm{d} \ell := -ik g^1_{k} - \partial_{y} g^2_{k}.
\end{align*}
For convenience, we denote $ \mathcal{A} =   - \nu (  \partial_y^2 -  k^2) + i k y  $.

For the periodic channel problem $\T \times [-1,1]$, for each fixed wave number $k$, the linearized operator $\mathcal{A}$ in the $y$-direction constitutes a Sturm--Liouville problem with a non-self-adjoint perturbation with a discrete spectrum. Due to the shear effect of the Couette flow, the enhanced dissipation rate is effectively amplified to $\sim \nu^{\frac13} |k|^{\frac23}$.

In the present paper, we consider the infinite channel $\R \times [-1,1]$, where the horizontal wave number $k$ varies continuously. The most critical difficulty arises as $k \to 0$ (corresponding to long-wave perturbations), for which the enhanced dissipation rate $\sim \nu^{\frac13} |k|^{\frac23}$ tends to zero. In this regime, the enhanced dissipation effect becomes extremely weak, and the dissipative time scale $\sim \nu^{-\frac13} |k|^{-\frac23}$ becomes very long. Likewise, the inviscid damping rate $\sim \frac{1}{|k| t^2}$ also slows down significantly. Consequently, nonlinear interactions can continually generate and excite these nearly undamped long-wave modes, making it difficult for the solution to decay toward equilibrium and potentially leading to instability.

\subsection{Space-Time Estimates for the Vorticity $ \omega_k$} Inspired by \cite{Miao2025Arxiv}, we decompose the horizontal wave numbers into three regimes: low (long-wave perturbations), intermediate, and high frequencies. For each range of the horizontal wave number $k$, we employ a combination of energy methods and Airy function analysis to derive relatively sharp space--time estimates corresponding to the three frequency regimes. More precisely, we have the following result:
\begin{proposition}\cite{Miao2025Arxiv}\label{w_resolvent_es}
	 Let $ \omega_k$ be a solution to system \eqref{equation_k} with initial data $ \omega_k^{\mathrm{in}}$, and assume that $ f^1_{k}, f^2_{k} \in L^2_y$.  Then there holds
\begin{eqnarray*}
E_k[\omega_k] \lesssim
	\left\{\begin{aligned}
		&\| \omega_k^{\mathrm{in}} \|_{L^2_y} + \nu^{-\frac{1}{6}} |k|^{\frac{2}{3}} \| (f^1_{k}, \theta_k)\|_{L^2_t L^2_y}+  \nu^{-\frac{1}{2}} \| f^2_{k}\|_{L^2_t L^2_y} ,\qquad  0 \le |k| \le 10 \nu ; \\
		&\| (\omega_k^{\mathrm{in}}, \nu^{\frac{1}{3} } |k|^{-\frac{1}{3}}\partial_{y} \omega_k^{\mathrm{in}} ) \|_{L^2_y} +  \nu^{-\frac{1}{6}} |k|^{\frac{2}{3}} \| (f^1_{k}, \theta_k)\|_{L^2_t L^2_y}\\
&\qquad\qquad\qquad\qquad\qquad\qquad\qquad+  \nu^{-\frac{1}{2}} \| f^2_{k}\|_{L^2_t L^2_y}, \qquad 10 \nu \le |k| \leq 1;  \\
		&\| (\omega_k^{\mathrm{in}},   |k|^{-1}\partial_{y} \omega_k^{\mathrm{in}} ) \|_{L^2_y} + \min\{\nu^{-\frac{1}{6}} |k|^{\frac{2}{3}}, \nu^{-\frac{1}{2}}\} \|(f^1_{k}, \theta_k)\|_{L^2_t L^2_y}\\
&\qquad\qquad\qquad\qquad\qquad\qquad\qquad\qquad\qquad+  \nu^{-\frac{1}{2}} \| f^2_{k}\|_{L^2_t L^2_y} ,  \qquad  |k| \geq 1.
	\end{aligned}\right.
\end{eqnarray*}

\end{proposition}
In addition, we consider the following non-self-adjoint perturbation of a Sturm--Liouville problem:
\begin{eqnarray}\label{eq_SL}
	\left\{\begin{aligned}
		&	- \nu(\partial_{y}^2 - k^2) f + ik(y- \lambda)f - \varepsilon \nu^{\frac{1}{3}} |k|^{\frac{2}{3}} f = F,\\
		&f|_{y = \pm 1} = 0.
	\end{aligned}\right.
\end{eqnarray}
We now state the following a priori resolvent estimates for the above system, mapping $ L^2 \rightarrow H^2$ and $ H^{-1} \rightarrow H^1$.
\begin{proposition}\label{prop_chen}
	 $ \bm{(L^2 \rightarrow H^2)}:$ Let $ f \in H^2_y$ be a solution of \eqref{eq_SL} with $ \lambda \in \R$ and $ F \in L^2_y$. Then, for any   $ k \in \R$,  there exists a constant $ C >0$ independent of $\nu$, $k$, and $ \lambda$, such that
	\begin{align*}
		\nu^{\frac{2}{3}} |k|^{\frac{1}{3}} \| \partial_{y} f\|_{L^2_y} + (\nu k^2)^{\frac{1}{3}} \|   f\|_{L^2_y} + |k| \| (y - \lambda) f\|_{L^2_y} \le C \|F \|_{L^2_y}.
	\end{align*}
	$ \bm{(H^{-1} \rightarrow H^1)}:$  Let $ f \in H^1_y$ be a solution of \eqref{eq_SL} with $ \lambda \in \R$ and $ F \in H^{-1}_y$. Then, for any   $ k \in \R$,  there exists a constant $ C >0$ independent of $\nu$, $k$, and $ \lambda$, such that
	\begin{align*}
		\nu \| \partial_{y} f \|_{L^2_y} + \nu^{\frac{2}{3}} |k|^{\frac{1}{3}} \|   f \|_{L^2_y} \le C \|F \|_{H^{-1}_y}.
	\end{align*}
	\begin{proof}
	We can establish these resolvent bounds by adapting the arguments developed in \cite{Chenqi2020ARMA}, specifically Proposition 3.1 and Lemma 3.4 therein.
	\end{proof}
\end{proposition}

\subsection{Space-time estimates for the temperature $ \theta_k$}
In this subsection, we decompose the second equation of system \eqref{equation_k} into homogeneous and nonhomogeneous components. Through this quasi-linearization approach, the ''bad'' nonlinear terms are isolated into an error system. By exploiting the zero initial data of the error system together with semigroup estimates for the linearized part, we can obtain certain decay properties, thus avoiding further restrictions on the size of the main system. More precisely, we decompose the temperature equation of $\theta_k$ into the homogeneous part $\theta^\mathrm{H}_k$ and inhomogeneous part $\theta^\mathrm{I}_k$:
\begin{eqnarray}\label{thetaH}
\left\{\begin{aligned}
&\partial_t \theta^\mathrm{H}_k + i k y \theta^\mathrm{H}_k  - \nu (  \partial_y^2 -  k^2) \theta^\mathrm{H}_k = 0,\\
&\theta^\mathrm{H}_k \big|_{t=0} =  \theta^{\mathrm{in}}_k(y), \quad \theta^\mathrm{H}_k \big|_{y = \pm 1} = 0,
\end{aligned}\right.
\end{eqnarray}
 and
 \begin{eqnarray}\label{eq_thetaL}
 \left\{\begin{aligned}
 &\partial_t \theta^\mathrm{I}_k + i k y \theta^\mathrm{I}_k  - \nu (  \partial_y^2 -  k^2) \theta^\mathrm{I}_k = -ik g^1_{k} - \partial_{y} g^2_{k},\\
 &\theta^\mathrm{I}_k \big|_{t=0} =  0, \quad \theta^\mathrm{I}_k \big|_{y = \pm 1} = 0.
 \end{aligned}\right.
 \end{eqnarray}

Next, we derive estimates for the homogeneous and nonhomogeneous components separately. We begin with the homogeneous part. Owing to the refined resolvent bounds and the enhanced dissipation structure of the linearized operator, we are able to obtain the following relatively sharp estimates.
\begin{lemma}\label{leta_h}
	Suppose that $ \theta_k^\mathrm{H}$ is a solution of system  \eqref{thetaH} with initial data $\theta^{\mathrm{in}}_k \in L^2_y $. Then, for any $ k \in \R$, the following estimate holds
	\begin{align*}
		\| \theta_k^\mathrm{H}\|_{L^2_y}  \le C e^{-c[(\nu k^2)^{\frac{1}{3}} + \nu] t} \| \theta^{\mathrm{in}}_k(y)\|_{L^2_y}.
	\end{align*}
	Moreover, it follows that
	\begin{align*}
	\big(\nu+ 	(\nu k^2)^{\frac{1}{3}}\big) \|  \theta_k^\mathrm{H} \|^2_{L^2_t L^2_y} \le C \| \theta^{\mathrm{in}}_k(y)\|^2_{L^2_y}.
	\end{align*}
	\begin{proof}
		First, we introduce the notation
		\begin{align*}
			(\mathcal{A} - i \lambda ) \theta_k^\mathrm{H} : = F
		\end{align*}
		for later convenience.
		
		Then, applying Proposition  \ref{prop_chen} to the homogeneous system, we immediately obtain
		\begin{align*}
			(\nu k^2)^{\frac{1}{3}} \|  \theta_k^\mathrm{H}  \|_{L^2_y} \lesssim \| F\|_{L^2_y}  =  \| (\mathcal{A} - i \lambda ) \theta_k^\mathrm{H}\|_{L^2_y}.
		\end{align*}
		According to the definition of $ \Phi(\mathcal{A})$(see Appendix), we obtain
		\begin{align}\label{A1}
			 \Phi(\mathcal{A}) \gtrsim (\nu k^2)^{\frac{1}{3}}.
		\end{align}
		On the other hand, by taking the $L^2-$inner product of $ \left(\mathcal{A} - i \lambda\right) \theta_k^\mathrm{H}$ with $ \theta_k^\mathrm{H}  $, we obtain
		\begin{align*}
			\operatorname{Re}\left\langle \left(\mathcal{A} - i \lambda\right) \theta_k^\mathrm{H}, \theta_k^\mathrm{H}   \right\rangle =  \nu k^2 \| \theta_k^\mathrm{H}\|^2_{L^2_y} + \nu \| \partial_{y} \theta_k^\mathrm{H}\|^2_{L^2_y} \le \| (\mathcal{A} - i \lambda ) \theta_k^\mathrm{H}\|_{L^2_y} \|   \theta_k^\mathrm{H}\|_{L^2_y}.
		\end{align*}
		Using the homogeneous boundary conditions satisfied by $ \theta_k^\mathrm{H}  $( $ \theta_k^\mathrm{H} \big|_{y = \pm 1} = 0$), and applying the Poincar\'e inequality, we naturally obtain
		\begin{align*}
			\| (\mathcal{A} - i \lambda ) \theta_k^\mathrm{H}\|_{L^2_y} \gtrsim \nu \| \theta_k^\mathrm{H}\|_{L^2_y}.
		\end{align*}
		Combining the above estimate with \eqref{A1}, we obtain
		\begin{align*}
			\Phi(\mathcal{A}) \gtrsim (\nu k^2)^{\frac{1}{3}} + \nu.
		\end{align*}
		It is straightforward to verify that $ \operatorname{Re} \langle \mathcal{A} f, f  \rangle \ge 0, $ for any $f$, which implies that $\mathcal{A}$ is an accretive operator. Moreover, for every $ \lambda < 0$, one can check that $ H = R(\lambda I - \mathcal{A})$, so that $\mathcal{A}$ is in fact $m-$accretive. Therefore, by invoking Lemma \ref{le-GP}, we conclude that
		\begin{align*}
		\| \theta_k^\mathrm{H}\|_{L^2_y} = \| e^{-\mathcal{A} t} \theta^{\mathrm{in}}_k(y) \|_{L^2_y} \le Ce^{\frac{\pi}{2}}  e^{-\Phi(\mathcal{A})} \| \theta^{\mathrm{in}}_k(y)\|_{L^2_y} \le C e^{-c[(\nu k^2)^{\frac{1}{3}} + \nu] t} \| \theta^{\mathrm{in}}_k(y)\|_{L^2_y}.
		\end{align*}
	This completes the proof of the first inequality in the lemma. The second inequality follows directly by integrating the above estimate in time.
	\end{proof}
\end{lemma}
Next, we proceed to derive the resolvent estimates for the inhomogeneous part of the temperature equation $ \theta_k^{\mathrm{I}}$. It is worth emphasizing that the underlying mechanisms differ substantially across frequency regimes; as a consequence, the inhomogeneous term exhibits markedly different behaviors in the high- and low-frequency ranges. Taking this distinction into account, we obtain the following estimates separately for the two regimes:
\begin{lemma}\label{le_thetaL}
	 Let $ \theta_k^\mathrm{I}$ be a solution of system \eqref{eq_thetaL}. For $ k$ in the medium-high frequency regime $ |k| \ge 10 \nu$, there exists a constant $C>0$, independent of $ \nu$, such that
	 \begin{align}\label{thetaI1_MH}
	 	\| \theta_k^\mathrm{I} \|_{L^\infty_t L^2_y} + 	(\nu k^2)^{\frac{1}{6}} \|  \theta_k^\mathrm{I} \|_{L^2_t L^2_y} \le C \big( \nu^{-\frac{1}{6}} |k|^{\frac{2}{3}} \| g^1_k\|_{L^2_t L^2_y} + \nu^{-\frac{1}{2}} \| g^2_k\|_{L^2_t L^2_y}\big).
	 \end{align}
	 In addition, for the low-frequency region $ 0 \le |k| \le 10 \nu$, one has
	 \begin{align}\label{thetaI1_L}
	 	\| \theta_k^\mathrm{I} \|_{L^\infty_t L^2_y} + 	\nu^{\frac{1}{2}} \|  \theta_k^\mathrm{I} \|_{L^2_t L^2_y} \le C \big( \nu^{-\frac{1}{6}} |k|^{\frac{2}{3}} \| g^1_k\|_{L^2_t L^2_y} + \nu^{-\frac{1}{2}} \| g^2_k\|_{L^2_t L^2_y}\big).
	 \end{align}
	 \begin{proof}
	 	To begin with, in the treatment of the medium- and high-frequency modes, we further refine the inhomogeneous contribution. Specifically, we decompose $ \theta_k^I$ into two components and examine separately the evolution equations satisfied by each part:
	 	\begin{eqnarray}\label{eq_thetaL12}
	 	\left\{\begin{aligned}
	 	&\partial_t \theta^{(\mathrm{I},1)}_k + i k y \theta^{(\mathrm{I},1)}_k  - \nu (  \partial_y^2 -  k^2) \theta^{(\mathrm{I},1)}_k = -ik g^1_{k}, \quad \theta^{(\mathrm{I},1)}_k \big|_{t=0} = \theta^{(\mathrm{I},1)}_k \big|_{y = \pm 1} = 0,\\
	 	& \partial_t \theta^{(\mathrm{I},2)}_k + i k y \theta^{(\mathrm{I},2)}_k  - \nu (  \partial_y^2 -  k^2) \theta^{(\mathrm{I},2)}_k =   - \partial_{y} g^2_{k}, \quad \theta^{(\mathrm{I},2)}_k \big|_{t=0} = \theta^{(\mathrm{I},2)}_k \big|_{y = \pm 1} = 0.
	 	\end{aligned}\right.
	 	\end{eqnarray}
	 	We then take the Fourier transform in the time variable $ t$. More precisely, we define $$  \Theta^{(\mathrm{I},i)}(\lambda, k, y) = \int_0^\infty \theta^{(\mathrm{I},i)}_k \cdot e^{-i\lambda t}  \mathrm{d} t, \quad \text{and} \quad  G^{i}(\lambda, k, y) = \int_0^\infty g^{i}_k \cdot e^{-i\lambda t}  \mathrm{d} t, \quad i = 1,2.$$
	 	Next, we apply the Fourier transform in the time variable $t$ to system \eqref{eq_thetaL12}. This yields the following transformed system:
	 	 \begin{eqnarray*}
	 	 \left\{\begin{aligned}
	 	 &  i k \left(y - \frac{\lambda}{k}\right) \Theta^{(\mathrm{I},1)}_k  - \nu (  \partial_y^2 -  k^2) \Theta^{(\mathrm{I},1)}_k = -ik G^1, \quad  \Theta^{(\mathrm{I},1)} \big|_{y = \pm 1} = 0,\\
	 	 &  i k \left(y - \frac{\lambda}{k}\right) \Theta^{(\mathrm{I},2)}_k  - \nu (  \partial_y^2 -  k^2) \Theta^{(\mathrm{I},2)}_k = -\partial_{y} G^2, \quad  \Theta^{(\mathrm{I},2)} \big|_{y = \pm 1} = 0.
	 	 \end{aligned}\right.
	 	 \end{eqnarray*}
	 	 By invoking Proposition \ref{prop_chen}, we obtain the corresponding estimates for both $ \Theta^{(\mathrm{I},1)}_k $ and $ \Theta^{(\mathrm{I},2)}_k $ separately.
	 	 \begin{align*}
	 	 	(\nu k^2)^{\frac{1}{3}} \|   \Theta^{(\mathrm{I},1)} \|_{L^2_y}   \le C \| kG^1 \|_{L^2_y} \quad \text{and} \quad \nu^{\frac{2}{3}} |k|^{\frac{1}{3}} \|   \Theta^{(\mathrm{I},2)} \|_{L^2_y} \le C \|G^2 \|_{L^2_y}.
	 	 \end{align*}
	 	 Next, by applying the Plancherel theorem together with the above inequalities, we obtain the following $L^2$-in-time estimate for $ \theta^{\mathrm{I}}_k$:
	 	 \begin{align}\label{first_thetaI}
	 	 		(\nu k^2)^{\frac{1}{6}} \|   \theta^{\mathrm{I}}_k \|_{L^2_t L^2_y} &\approx 	(\nu k^2)^{\frac{1}{6}} \big(\|   \Theta^{(\mathrm{I},1)} \|_{L^2_\lambda L^2_y} + \|   \Theta^{(\mathrm{I},2)} \|_{L^2_\lambda L^2_y} \big)\nn\\
	 	 		& \lesssim \nu^{-\frac{1}{6}} |k|^{\frac{2}{3}} \| G^1\|_{L^2_\lambda L^2_y} + \nu^{-\frac{1}{2}} \| G^2\|_{L^2_\lambda L^2_y}\nn\\
	 	 		& \approx \nu^{-\frac{1}{6}} |k|^{\frac{2}{3}} \| g^1_k\|_{L^2_t L^2_y} + \nu^{-\frac{1}{2}} \| g^2_k\|_{L^2_t L^2_y}.
	 	 \end{align}
	 	 For the $L^\infty$-in-time estimate of $\theta_k^\mathrm{I}$, we can directly perform an $L^2$ energy estimate on system \eqref{eq_thetaL}, which yields:
	 	 \begin{align*}
	 	 	\operatorname{Re} \left\langle \left( \partial_t - \mathcal{A}\right) \theta^{\mathrm{I}}_k,  \theta^{\mathrm{I}}_k \right\rangle = - \operatorname{Re} \left\langle  ik g^1_k + \partial_{y} g^2_k,  \theta^{\mathrm{I}}_k \right\rangle.
	 	 \end{align*}
	 	 The left-hand side of the above expression can be written as $$ \text{LHS} = \frac{1}{2} \frac{d}{dt} \| \theta^{\mathrm{I}}_k\|^2_{L^2_y} + \nu \| \partial_{y} \theta^{\mathrm{I}}_k\|^2_{L^2_y} + \nu k^2 \| \theta^{\mathrm{I}}_k\|^2_{L^2_y}.$$
	 	The right-hand side of the above expression can be controlled by applying Young's inequality and integration by parts, yielding $$\text{RHS}  \le   \frac{1}{4} (\nu k^2)^{\frac{1}{3}} \| \theta^{\mathrm{I}}_k\|^2_{L^2_y} + \nu^{-\frac{1}{3}} |k|^{\frac{4}{3}} \| g^1_k\|^2_{L^2_y} + \frac{1}{4} \nu \| \partial_{y} \theta^{\mathrm{I}}_k\|^2_{L^2_y} +  \nu^{-1} \| g^2_k\|^2_{L^2_y}.$$
	 	 Combining the result of \eqref{first_thetaI}, we immediately obtain the estimate in \eqref{thetaI1_MH}. In fact, this estimate holds for all  $ k \in \R$.
	 	
	 	 Next, we require stronger control for the low-frequency regime. When $k=0$, \eqref{thetaI1_L} clearly holds.
	 	 For $k \neq 0$, a direct $L^2$ energy estimate gives
	 	 \begin{align*}
	 	 	\frac{1}{2} \frac{d}{dt} \| \theta^{\mathrm{I}}_k\|^2_{L^2_y} + \frac{3}{4}\nu \| \partial_{y} \theta^{\mathrm{I}}_k\|^2_{L^2_y} + \nu k^2 \| \theta^{\mathrm{I}}_k\|^2_{L^2_y} \le \frac{1}{4} (\nu k^2)^{\frac{1}{3}} \| \theta^{\mathrm{I}}_k\|^2_{L^2_y} + \nu^{-\frac{1}{3}} |k|^{\frac{4}{3}} \| g^1_k\|^2_{L^2_y}   +  \nu^{-1} \| g^2_k\|^2_{L^2_y}.
	 	 \end{align*}
	 	 Using the homogeneous boundary conditions satisfied by $ \theta^{\mathrm{I}}_k$, we can apply the Poincar\'e inequality
	 	 $$ \| \partial_{y} \theta^{\mathrm{I}}_k\|^2_{L^2_y} \gtrsim \|   \theta^{\mathrm{I}}_k\|^2_{L^2_y},$$ which completes the proof of the lemma.
	 \end{proof}
\end{lemma}

Combining the results of Lemma \ref{leta_h} and Lemma \ref{le_thetaL}, we obtain the resolvent estimates for $\theta_k$ in both the high- and low-frequency regimes. More precisely, we have the following lemma.
\begin{lemma}\label{le_theta}
		Let $ \theta_k$ be a solution of system \eqref{equation_k}. For $ k$ in the medium-high frequency regime $ |k| \ge 10 \nu$, there exists a constant $C>0$, independent of $ \nu$, such that
		\begin{align*}
		\| \theta_k  \|_{L^\infty_t L^2_y} + 	(\nu k^2)^{\frac{1}{6}} \|  \theta_k  \|_{L^2_t L^2_y} \le C\| \theta_k^{\mathrm{in}} \|_{L^2_y}  + C \big( \nu^{-\frac{1}{6}} |k|^{\frac{2}{3}} \| g^1_k\|_{L^2_t L^2_y} + \nu^{-\frac{1}{2}} \| g^2_k\|_{L^2_t L^2_y}\big).
		\end{align*}
		In addition, for the low-frequency region $ 0 \le |k| \le 10 \nu$, one has
		\begin{align*}
		\| \theta_k  \|_{L^\infty_t L^2_y} + 	\nu^{\frac{1}{2}} \|  \theta_k  \|_{L^2_t L^2_y} \le C\| \theta_k^{\mathrm{in}} \|_{L^2_y}  + C \big( \nu^{-\frac{1}{6}} |k|^{\frac{2}{3}} \| g^1_k\|_{L^2_t L^2_y} + \nu^{-\frac{1}{2}} \| g^2_k\|_{L^2_t L^2_y}\big).
		\end{align*}
\end{lemma}

\section{Nonlinear stability analysis for the Boussinesq equations}
In this section, our primary goal is to control the nonlinear terms and the buoyancy contribution in the velocity equation. To this end, we begin by decomposing the frequencies in the $x-$direction into three ranges: low, medium, and high frequencies. Since the nonlinear interactions involve products of modes from different frequency ranges, the decomposition leads to nine possible interaction scenarios. This separation enables us to take advantage of the distinct dissipation or mixing mechanisms available in each frequency regime, thereby allowing for a more effective control of the nonlinear contributions.

More precisely, we partition the horizontal frequencies into the following three regions:
$$ \R = \cup_{i=1}^{3} I_i,\quad  I_1 = \big\{ k \big| |k| \le 10 \nu \big\}, \quad I_2 =\left\{ k \big| 10 \nu \le |k| \le 1  \right\}, \quad I_3 = \left\{ k \big| |k| \ge 1 \right\}. $$
To further control the buoyancy term in the vorticity equation, it is necessary to refine the analysis within the high-frequency regime $I_3$. In particular, we decompose this region into the following two subranges:
$$I_3 =   I_4 \cup I_5, \quad\hbox{with}\quad  I_4 = \left\{ k \big|1 \le  |k| \le  \nu^{-\frac{1}{2}} \right\} \quad\hbox{and}\quad I_5 = \left\{ k \big|   |k| \ge \nu^{-\frac{1}{2}} \right\}.$$
Based on the definitions of  $ E_k[\omega_k]$ and $ E_k[\theta_k]$ in the respective frequency regions, we can apply the Gagliardo--Nirenberg inequality to obtain
\begin{align}\label{u1_es}
&\left\| u^1_k \right\|_{L^\infty_t L^{\infty}_y} \lesssim \min\{ 1, |k|^{-\frac{1}{2}}\}  E_k[\omega_k], \quad\left\|u^1_k\right\|_{L^2_t L^\infty_y} \lesssim \begin{cases}\nu^{-\frac{1}{2}} E_k[\omega_k], & k \in I_1; \\ \nu^{-\frac{1}{8}}|k|^{-\frac{3}{8}} E_k[\omega_k], & k \in I_2; \\ \nu^{-\frac{1}{8}}|k|^{-\frac{3}{4}} E_k[\omega_k], & k \in I_3;\end{cases}
\end{align}
and
\begin{align}\label{theta_es}
&\left\| \theta_k \right\|_{L^\infty_t L^{2}_y} \lesssim  \begin{cases}  E_k[\theta_k], & k \in I_1\cup I_2; \\  |k|^{-\frac{1}{3}} E_k[\theta_k], & k \in I_3; \end{cases}, \quad
\left\|\theta_k\right\|_{L^2_t L^2_y} \lesssim \begin{cases}\nu^{-\frac{1}{2}} E_k[\theta_k], & k \in I_1; \\ \nu^{-\frac{1}{6}}|k|^{-\frac{1}{3}} E_k[\theta_k], & k \in I_2; \\ \nu^{-\frac{1}{6}}|k|^{-\frac{2}{3}} E_k[\theta_k], & k \in I_3.\end{cases}
\end{align}
For the convenience of the subsequent analysis, we first introduce the following two energy functionals: $$ \mathcal{E}_{\omega} = 	\big\| E_k[\omega_k] \big\|_{L^1_k(\R)} \quad \text{and} \quad \mathcal{E}_{\theta} = 	\big\| E_k[\theta_k] \big\|_{L^1_k(\R)}.$$
Based on the introduced energy functionals and by applying Proposition \ref{w_resolvent_es}, the nonlinear terms $ f^1_k$ and $ f^2_k$ can be controlled separately following the approach developed in \cite{Miao2025Arxiv}. Consequently, we obtain
 \begin{align*}
 	\big\| E_k[\omega_k] \big\|_{L^1_k(\R)}
 &\lesssim \left\| E_k^{\mathrm{in}}[\omega_k] \right\|_{L^1_k(\R)} + \left\| \min\{ \nu^{-\frac{1}{6}} |k|^{\frac{2}{3}} , \nu^{-\frac{1}{2}} \} \left(f^1_k, \theta_k\right)\right\|_{L^1_k L^2_t L^2_y} + \nu^{-\frac{1}{2}} \left\|     f^2_k\right\|_{L^1_k L^2_t L^2_y}\\
 	&\lesssim \left\| E_k^{\mathrm{in}}[\omega_k] \right\|_{L^1_k(\R)} +  \nu^{-\frac{1}{2}} \big\| E_k[\omega_k] \big\|^2_{L^1_k(\R)} +  \left\| \min\{ \nu^{-\frac{1}{6}} |k|^{\frac{2}{3}} , \nu^{-\frac{1}{2}} \}    \theta_k \right\|_{L^1_k L^2_t L^2_y}.
 \end{align*}
In the treatment of the buoyancy term described above, by combining \eqref{theta_es} with the frequency decomposition applied to $I_3$, we obtain
\begin{align*}
	\left\| \min\{ \nu^{-\frac{1}{6}} |k|^{\frac{2}{3}} , \nu^{-\frac{1}{2}} \}    \theta_k \right\|_{L^1_k L^2_t L^2_y} &\lesssim  \begin{cases}\nu^{-\frac{1}{6}} \nu^{-\frac{1}{2}} \left\| |k|^{\frac{2}{3}} E_k[\theta_k]\right\|_{L^1_k(I_1)}, & k \in I_1; \\ \nu^{-\frac{1}{6}} \nu^{-\frac{1}{6}} \left\| |k|^{\frac{1}{3}} E_k[\theta_k]\right\|_{L^1_k(I_2)}, & k \in I_2; \\ \nu^{-\frac{1}{6}} \nu^{-\frac{1}{6}}  \left\|   E_k[\theta_k]\right\|_{L^1_{k }(I_4)}, & k \in I_4;\\ \nu^{-\frac{1}{6}} \nu^{-\frac{1}{2}}  \left\| |k|^{-\frac{2}{3}}  E_k[\theta_k]\right\|_{L^1_{k }(I_5)}, & k \in I_5;\end{cases}\\
	&\lesssim \nu^{-\frac{1}{3}} 	\big\| E_k[\theta_k] \big\|_{L^1_k(\R)}.
\end{align*}
As a result, we obtain the following energy estimate, in which the threshold for the temperature needs only to be lower than that of the velocity by $ O(\nu^{\frac13})$:
\begin{align}\label{omega_es}
\big\| E_k[\omega_k] \big\|_{L^1_k(\R)} \lesssim \left\| E_k^{\mathrm{in}}[\omega_k] \right\|_{L^1_k(\R)} +  \nu^{-\frac{1}{2}} \big\| E_k[\omega_k] \big\|^2_{L^1_k(\R)} + \nu^{-\frac{1}{3}} 	\big\| E_k[\theta_k] \big\|_{L^1_k(\R)}.
\end{align}
Therefore, it remains to estimate the nonlinear terms in the temperature equation. According to Proposition \ref{le_theta}, the behavior of $\theta_k$ differs across the low-, medium-, and high-frequency regimes.
Accordingly, we divide the nonlinear estimates for the temperature into three subparts--corresponding to the low-, medium-, and high-frequency ranges and analyze each part separately.
\subsection{Nonlinear estimates for low frequency}
For the low-frequency regime, specifically, we have
\begin{lemma}\label{prop_theta_low}
	Let $ \theta_k$ be a solution to system \eqref{equation_k}.
	For the low-frequency regime $ |k| \le 10 \nu$, there exists a constant $C$, independent of $\nu$ and $k$, such that
	\begin{align*}
	\left\|  \| \theta_k  \|_{L^\infty_t L^2_y} \right\|_{L^1_{|k| \le 10 \nu}} + 	\nu^{\frac{1}{2}} \left\| \|  \theta_k  \| _{L^2_t L^2_y} \right\|_{ |k| \le 10 \nu}  \le C\left\| \left\| \theta_k^{\mathrm{in}} \right\|_{L^2_y} \right\|_{L^1_{|k| \le 10 \nu}} + C \nu^{-\frac{1}{2}} \mathcal{E}_{\omega} \cdot \mathcal{E}_{\theta}.
	\end{align*}
\begin{proof}
	Based on the definition of $ E_k[\theta_k ]$	and in combination with Lemma \ref{le_theta}, we can directly perform an energy estimate to obtain
	\begin{align}\label{ta1}
			\big\| E_k[\theta_k ] \big\|_{L^1_{|k| \le 10 \nu}} &\lesssim \left\| \left\| \theta_k^{\mathrm{in}} \right\|_{L^2_y} \right\|_{L^1_{|k| \le 10 \nu}} +  \int_{|k| \le 10 \nu}  \left( \nu^{-\frac{1}{6}} |k|^{\frac{2}{3}} \left\| g^1_k \right\|_{L^2_t L^2_y} +  \nu^{-\frac{1}{2}}  \left\| g^2_k \right\|_{L^2_t L^2_y}  \right) \mathrm{d} k\nn\\
		&\lesssim \left\| \left\| \theta_k^{\mathrm{in}} \right\|_{L^2_y} \right\|_{L^1_{|k| \le 10 \nu}} + R_1^L + R_2^L.
	\end{align}
	Next, we estimate the terms $  R_1^L$	and $  R_2^L$	individually. For $ R_2^L$, by applying the Gagliardo--Nirenberg inequality $ \left\|u^2_\ell   \right\|_{L^2_t L^\infty_y} \lesssim \left\|u^2_\ell   \right\|^{\frac{1}{2}}_{L^2_t L^2_y} \cdot \left\|\partial_{y} u^2_\ell   \right\|^{\frac{1}{2}}_{L^2_t L^2_y}  $ together with the Poincar\'e inequality $  \left\|u^2_\ell   \right\|_{L^2_t L^2_y} \lesssim \left\|\partial_{y} u^2_\ell   \right\|_{L^2_t L^2_y}$, it can be decomposed into three components as follows:
	\begin{align*}
		R_2^L  \lesssim& \nu^{-\frac{1}{2}} \sum_{i = 1}^{3}  \int_{I_1}  \int_{I_i} \left\|u^2_\ell   \right\|_{L^2_t L^\infty_y}  \left\|  \theta_{k-\ell} \right\|_{L^\infty_t L^2_y}  \mathrm{d} \ell  \mathrm{d} k\\
		\lesssim& \nu^{-\frac{1}{2}} \sum_{i = 1}^{3}  \int_{I_1}  \int_{I_i}  \left\|\partial_{y} u^2_\ell   \right\|_{L^2_t L^2_y}   \cdot \left\|  \theta_{k-\ell} \right\|_{L^\infty_t L^2_y}  \mathrm{d} \ell  \mathrm{d} k\\
		\lesssim& \nu^{-\frac{1}{2}} \sum_{i = 1}^{3}  \int_{I_1}  \int_{I_i}  |\ell|\left\|  u^1_\ell   \right\|_{L^2_t L^2_y}   \cdot \left\|  \theta_{k-\ell} \right\|_{L^\infty_t L^2_y}  \mathrm{d} \ell  \mathrm{d} k\\
 :=& R_{2,1}^L + R_{2,2}^L + R_{2,3}^L.
	\end{align*}
	For the first component $ R_{2,1}^L$, using \eqref{theta_es} and the fact that  $ \left\| \ell u^1_\ell \right\|_{L^2_tL^2_y} \lesssim E_k[\omega_k]$, we obtain
	\begin{align*}
		R_{2,1}^L \lesssim \nu^{-\frac{1}{2}}  \int_{I_1}  \int_{I_1}  E_\ell[\omega_\ell]      \cdot  E_{k-\ell}[\theta_{k-\ell}]  \mathrm{d} \ell  \mathrm{d} k \lesssim \nu^{-\frac{1}{2}} \mathcal{E}_{\omega} \cdot \mathcal{E}_{\theta}.
	\end{align*}
	For the second component $ R_{2,2}^L$, in view of Lemma \ref{lebasis1}, it is straightforward to deduce
	\begin{align*}
		 R_{2,2}^L &\lesssim \nu^{-\frac{1}{2}}  \left(\int_{0 \le |k-\ell| \le 10 \nu} + \int_{10 \nu \le |k-\ell| \le 2} \right) \int_{10 \nu \le |\ell| \le 1} |\ell|\left\|  u^1_\ell   \right\|_{L^2_t L^2_y}   \cdot \left\|  \theta_{k-\ell} \right\|_{L^\infty_t L^2_y}   \mathrm{d} \ell  \mathrm{d} k\\
		 &\lesssim \nu^{-\frac{1}{2}} \mathcal{E}_{\omega} \cdot \mathcal{E}_{\theta}.
	\end{align*}
	For the third component $ R_{2,3}^L$, since $ |k| \le 10 \nu$ and $ |\ell| \ge 1$, it can be further decomposed into three subregions for separate treatment.
	Subsequently, by applying \eqref{theta_es} together with the fact that $ \left\| \ell u^1_\ell \right\|_{L^2_tL^2_y} \lesssim E_k[\omega_k]$, we obtain
	\begin{align*}
		R_{2,3}^L &\lesssim \nu^{-\frac{1}{2}}  \left(\int_{0 \le |k-\ell| \le 10 \nu} + \int_{10 \nu \le |k-\ell| \le 1}  + \int_{  |k-\ell| \ge 1}\right) \int_{  |\ell| \ge 1} |\ell|\left\|  u^1_\ell   \right\|_{L^2_t L^2_y}   \cdot \left\|  \theta_{k-\ell} \right\|_{L^\infty_t L^2_y}   \mathrm{d} \ell  \mathrm{d} k\\
	  &\lesssim \nu^{-\frac{1}{2}} \mathcal{E}_{\omega} \cdot \mathcal{E}_{\theta}.
	\end{align*}
	
	Next, we turn to the treatment of $ R_1^L$, Similar to the approach for
	$R_2^L $, we decompose $ R_1^L$ into three subparts according to frequency:
	\begin{align*}
		R_1^L \lesssim \nu^{-\frac{1}{6}} \sum_{i = 1}^{3} \int_{  0 \le |k| \le 10 \nu  } \int_{I_i}  |k|^{\frac{2}{3}} \left\| u^1_\ell \right\|_{L^2_t L^\infty_y}  \left\| \theta_{k-\ell} \right\|_{L^\infty_t L^2_y}  \mathrm{d} \ell \mathrm{d}k := 	R_{1,1}^L + 	R_{1,2}^L + 	R_{1,3}^L.
	\end{align*}
	For $R_{1,1}^L $, when $  0 \le |k| \le 10 \nu$ and $  0 \le |\ell| \le 10 \nu$, we have $  0 \le |k-\ell| \lesssim  10 \nu$. By combining \eqref{u1_es}, \eqref{theta_es}, and H\"older's inequality, we obtain	
	\begin{align*}
		R_{1,1}^L &\lesssim \nu^{-\frac{1}{6}}   \int_{  0 \le |k-\ell| \lesssim  10 \nu  } \int_{I_1}  |\nu|^{\frac{2}{3}} \left\| u^1_\ell \right\|_{L^2_t L^\infty_y}  \left\| \theta_{k-\ell} \right\|_{L^\infty_t L^2_y}  \mathrm{d} \ell \mathrm{d}k\\
		&\lesssim \nu^{-\frac{1}{2}} \mathcal{E}_{\omega} \cdot \mathcal{E}_{\theta}.
	\end{align*}
	For $R_{1,2}^L $, when $  0 \le |k| \le 10 \nu$ and $  10 \nu \le |\ell| \le 1$, we have $  0 \le |k-\ell| \lesssim  1$. Using \eqref{u1_es} and the fact that $ |k| \cdot |\ell|^{-1} \le 1$, we deduce
	\begin{align*}
		R_{1,2}^L &\lesssim \nu^{-\frac{1}{6} - \frac{1}{8}}   \int_{  0 \le |k-\ell| \lesssim  1   } \int_{10 \nu \le |\ell| \le 1}  |k|^{\frac{2}{3}} |\ell|^{-\frac{3}{8}} E_\ell[\omega_\ell]      \cdot  E_{k-\ell}[\theta_{k-\ell}]   \mathrm{d} \ell \mathrm{d}k\\
		&\lesssim  \nu^{-\frac{7}{24}} \mathcal{E}_{\omega} \cdot \mathcal{E}_{\theta}.
	\end{align*}
	For $ 	R_{1,3}^L $, when $  0 \le |k| \le 10 \nu$ and $    |\ell| \ge 1$, noting that $ 10 \nu  \le 10 \nu_0 \le \frac{1}{2}$, we have $   |k-\ell| \approx |\ell|$. Applying \eqref{u1_es} once more, we finally obtain
	\begin{align*}
		R_{1,3}^L &\lesssim \nu^{-\frac{1}{6} -\frac{1}{8}}   \int_{    |k-\ell| \gtrsim  1   } \int_{ |\ell| \ge 1}  |k|^{\frac{2}{3}} |\ell|^{-\frac{3}{4}} E_\ell[\omega_\ell]      \cdot  E_{k-\ell}[\theta_{k-\ell}]   \mathrm{d} \ell \mathrm{d}k\\
		&\lesssim  \nu^{-\frac{7}{24}} \mathcal{E}_{\omega} \cdot \mathcal{E}_{\theta}.
	\end{align*}
	Combining all the above estimates, we complete the proof of the lemma.
\end{proof}
	
\end{lemma}

\subsection{Nonlinear estimates for medium frequency}
For the medium-frequency regime, specifically, we have
\begin{lemma}\label{prop_theta_medium}
	Let $ \theta_k$ be a solution to system \eqref{equation_k}.
	For the medium-frequency regime $ 10 \nu \le  |k| \le 1$, there exists a constant $C$, independent of $\nu$ and $k$, such that
	\begin{align*}
	\left\|  \| \theta_k  \|_{L^\infty_t L^2_y} \right\|_{L^1_{10 \nu \le |k| \le 1}} + 	\nu^{\frac{1}{6}} \left\| \big\| |k|^{\frac{1}{3}} \theta_k  \big\| _{L^2_t L^2_y} \right\|_{10 \nu \le |k| \le 1 }  \le C\left\| \left\| \theta_k^{\mathrm{in}} \right\|_{L^2_y} \right\|_{L^1_{10 \nu \le |k| \le 1 }} + C \nu^{-\frac{1}{2}} \mathcal{E}_{\omega} \cdot \mathcal{E}_{\theta}.
	\end{align*}
	\begin{proof}
		By applying an estimate analogous to that of \eqref{ta1}, we obtain
		\begin{align*}
		\big\| E_k[\theta_k ] \big\|_{L^1_{10 \nu \le |k| \le 1}} &\lesssim \left\| \left\| \theta_k^{\mathrm{in}} \right\|_{L^2_y} \right\|_{L^1_{10 \nu \le |k| \le 1}} +  \int_{10 \nu \le |k| \le 1}  \left( \nu^{-\frac{1}{6}} |k|^{\frac{2}{3}} \left\| g^1_k \right\|_{L^2_t L^2_y} +  \nu^{-\frac{1}{2}}  \left\| g^2_k \right\|_{L^2_t L^2_y}  \right) \mathrm{d} k\\
		&\lesssim \left\| \left\| \theta_k^{\mathrm{in}} \right\|_{L^2_y} \right\|_{L^1_{10 \nu \le |k| \le 1 }} + R_1^M + R_2^M.
		\end{align*}
		Next, we turn to the treatment of $ R_2^M$.
		Following the same strategy as in the analysis of $ R_2^L$, we decompose $R_2^M $	into three parts:
		\begin{align}\label{R2_M}
		R_2^M &\lesssim \nu^{-\frac{1}{2}} \int_{10 \nu \le |k| \le 1}    \left\| u^2_\ell \cdot \theta_{k-\ell} \right\|_{L^2_t L^2_y}    \mathrm{d} \ell   \mathrm{d} k \nn\\
		&\lesssim \nu^{-\frac{1}{2}} \sum_{i = 1}^{3}  \int_{I_2}  \int_{I_i}    \left\|  u^2_\ell   \right\|^{\frac{1}{2}}_{L^2_t L^2_y} \left(|\ell|\left\|  u^1_\ell   \right\|_{L^2_t L^2_y}\right)^{\frac{1}{2}}   \cdot \left\|  \theta_{k-\ell} \right\|_{L^\infty_t L^2_y}  \mathrm{d} \ell  \mathrm{d} k := R_{2,1}^M + R_{2,2}^M + R_{2,3}^M.
		\end{align}
		For $ R_{2,1}^M$	and $ R_{2,2}^M$, when  $ 10 \nu \le |k| \le 1$ and $ 0 \le |\ell| \le 1 $, we have $ 0 \le |k-\ell| \lesssim 1$. Using estimate \eqref{theta_es} together with $ \| k \u_k \|_{L^2_tL^2_y} \lesssim E_k[\omega_k]$, we readily obtain
		\begin{align*}
			R_{2,1}^M + R_{2,2}^M \lesssim \nu^{-\frac{1}{2}} \int_{0 \le |k-\ell| \lesssim  1} \int_{  0 \le |\ell| \le 1  } |\ell| \left\| \u_\ell \right\|_{L^2_t L^2_y}     \left\|  \theta_{k-\ell} \right\|_{L^\infty_t L^2_y}      \mathrm{d} \ell  \mathrm{d} k \lesssim \nu^{-\frac{1}{2}} \mathcal{E}_{\omega} \cdot \mathcal{E}_{\theta}.
		\end{align*}
		For $ R_{2,3}^M$, using \eqref{theta_es} directly yields
		\begin{align*}
			R_{2,3}^M  &\lesssim \nu^{-\frac{1}{2}} \left( \int_{0 \le |k-\ell| \le   1} + \int_{  |k-\ell| \ge  1}\right) \int_{  |\ell| \ge 1  } |\ell| \left\| \u_\ell \right\|_{L^2_t L^2_y}     \left\|  \theta_{k-\ell} \right\|_{L^\infty_t L^2_y}      \mathrm{d} \ell  \mathrm{d} k\\
			&\lesssim \nu^{-\frac{1}{2}} \mathcal{E}_{\omega} \cdot \mathcal{E}_{\theta}.
		\end{align*}
		Next, we turn to the treatment of the more delicate term $ R_1^M$. We first perform a frequency decomposition and rewrite $ R_1^M$	as the sum of the following three components:
		\begin{align*}
		R_1^M \lesssim \nu^{-\frac{1}{6}} \sum_{i = 1}^{3} \int_{ 10 \nu \le |k| \le  1 } \int_{I_i}  |k|^{\frac{2}{3}} \left\| u^1_\ell \cdot \theta_{k-\ell}\right\|_{L^2_t L^2_y}    \mathrm{d} \ell \mathrm{d}k := 	R_{1,1}^M + 	R_{1,2}^M + 	R_{1,3}^M.
		\end{align*}
		For the term $ 	R_{1,1}^M$, when  $ 10 \nu \le |k| \le  1$ and $ 0 \le |\ell| \le  10 \nu$. In this regime we have $0 \le  |k-\ell| \lesssim 1$, so that $ u^1_\ell$	lies in the low--frequency region, while $ k-\ell$ remains in the mid-low--frequency regime.
		Consequently, we may directly invoke the estimates in \eqref{u1_es} and \eqref{theta_es} to obtain
		\begin{align*}
			R_{1,1}^M &\lesssim \nu^{-\frac{1}{6}} \left(\int_{0 \le |k-\ell| \le   10 \nu} + \int_{10 \nu \le |k-\ell| \lesssim  1}\right) \int_{  0 \le |\ell| \le 10 \nu  } |k|^{\frac{2}{3}} \left\| u^1_\ell \cdot \theta_{k-\ell}\right\|_{L^2_t L^2_y}      \mathrm{d} \ell  \mathrm{d} k \\
			&\lesssim \nu^{-\frac{1}{6}}  \int_{0 \le |k-\ell| \le   10 \nu}   \int_{  0 \le |\ell| \le 10 \nu  } \left(|k-\ell|^{\frac{2}{3}} + |\ell|^{\frac{2}{3}} \right) \left\| u^1_\ell \right\|_{L^2_t L^\infty_y}     \left\|  \theta_{k-\ell} \right\|_{L^\infty_t L^2_y}      \mathrm{d} \ell  \mathrm{d} k \\
			&\quad + \nu^{-\frac{1}{6}}  \int_{10 \nu \le |k-\ell| \lesssim  1}  \int_{  0 \le |\ell| \le 10 \nu  } \left(|k-\ell|^{\frac{2}{3}} + |\ell|^{\frac{2}{3}} \right)  \left\| u^1_\ell \right\|_{L^\infty_t L^\infty_y}     \left\|  \theta_{k-\ell} \right\|_{L^2_t L^2_y}      \mathrm{d} \ell  \mathrm{d} k \\
			&\lesssim \mathcal{E}_{\omega} \cdot \mathcal{E}_{\theta} + \nu^{-\frac{1}{3}}  \int_{10 \nu \le |k-\ell| \lesssim  1}  \int_{  0 \le |\ell| \le 10 \nu  } \left(|k-\ell|^{\frac{1}{3}} + \frac{|\ell|^{\frac{2}{3}}}{|k-\ell|^{\frac{1}{3}}} \right)   E_\ell[\omega_\ell]      \cdot  E_{k-\ell}[\theta_{k-\ell}]       \mathrm{d} \ell  \mathrm{d} k \\
			&\lesssim \nu^{-\frac{1}{3}} \mathcal{E}_{\omega} \cdot \mathcal{E}_{\theta}.
		\end{align*}
		To handle $ R_{1,2}^M$, we perform a further decomposition, which is given by the following partition.
		\begin{align*}
			R_{1,2}^M  \lesssim& \nu^{-\frac{1}{6}} \left(\int_{0 \le |k-\ell| \le   10 \nu} + \int_{10 \nu \le |k-\ell| \lesssim  1}\right) \int_{  10 \nu \le |\ell| \le  1 } |k|^{\frac{2}{3}} \left\| u^1_\ell \cdot \theta_{k-\ell}\right\|_{L^2_t L^2_y}      \mathrm{d} \ell  \mathrm{d} k \\
:=&  R_{1,2,1}^M  + R_{1,2,2}^M.
		\end{align*}
		The subterm $ R_{1,2,1}^M$	can be treated in exactly the same manner as $ R_{1,1}^M$.
		Note that in this region the inequality $ \frac{|k-\ell|^{\frac{2}{3}} + |\ell|^{\frac{2}{3}}}{|\ell|^{\frac{3}{8}}} \lesssim 1$ holds. Therefore, by combining estimates \eqref{u1_es} and \eqref{theta_es}, we obtain
		\begin{align*}
			R_{1,2,1}^M \lesssim& \nu^{-\frac{1}{6}}  \int_{0 \le |k-\ell| \le   10 \nu}   \int_{   10 \nu \le |\ell| \le  1  } \left(|k-\ell|^{\frac{2}{3}} + |\ell|^{\frac{2}{3}} \right) \left\| u^1_\ell \right\|_{L^2_t L^\infty_y}     \left\|  \theta_{k-\ell} \right\|_{L^\infty_t L^2_y}      \mathrm{d} \ell  \mathrm{d} k  \\
\lesssim& \nu^{-\frac{7}{24}} \mathcal{E}_{\omega} \cdot \mathcal{E}_{\theta}.
		\end{align*}
		For $ R_{1,2,2}^M$, we apply the inequality $ |k| \le |k-\ell| + |\ell|$ to further decompose it into two terms:
		\begin{align*}
			R_{1,2,2}^M \lesssim& \nu^{-\frac{1}{6}}   \int_{10 \nu \le |k-\ell| \lesssim  1}  \int_{  10 \nu \le |\ell| \le  1 } \left(|k-\ell|^{\frac{2}{3}} + |\ell|^{\frac{2}{3}} \right)  \left\| u^1_\ell \cdot \theta_{k-\ell}\right\|_{L^2_t L^2_y}      \mathrm{d} \ell  \mathrm{d} k \\
: =& R_{1,2,2,1}^M   +R_{1,2,2,2}^M.
		\end{align*}
		The term $ R_{1,2,2,1}^M $	can be directly estimated using the energy estimates \eqref{u1_es} and \eqref{theta_es}
		\begin{align*}
			R_{1,2,2,1}^M \lesssim& \nu^{-\frac{1}{6}}   \int_{10 \nu \le |k-\ell| \lesssim  1}  \int_{  10 \nu \le |\ell| \le  1 }  |k-\ell|^{\frac{2}{3}}  \left\| u^1_\ell \right\|_{L^\infty_t L^\infty_y}     \left\|  \theta_{k-\ell} \right\|_{L^2_t L^2_y}        \mathrm{d} \ell  \mathrm{d} k \\
\lesssim& \nu^{-\frac{1}{3}} \mathcal{E}_{\omega} \cdot \mathcal{E}_{\theta}.
		\end{align*}
		 When $ 10 \nu \le |k-\ell| \lesssim  1$ and $ 10 \nu \le |\ell| \le  1 $, we have $ |k-\ell|^{-\frac{1}{6}} \lesssim \nu^{-\frac{1}{6}}$, combining estimates \eqref{u1_es} and \eqref{theta_es}, we obtain
		 \begin{align}\label{ta2}
		 \nu^{-\frac{1}{6}} \left\| u^1_\ell \cdot \theta_{k-\ell}\right\|_{L^2_t L^2_y} &  \lesssim \nu^{-\frac{1}{6}}\left( \left\| u^1_\ell \right\|_{L^\infty_t L^\infty_y}     \left\|  \theta_{k-\ell} \right\|_{L^2_t L^2_y}  \right)^{\frac{5}{9}}  \cdot \left(  \left\| u^1_\ell \right\|_{L^2_t L^\infty_y}     \left\|  \theta_{k-\ell} \right\|_{L^\infty_t L^2_y}    \right)^{\frac{4}{9}}\nn\\
		 &\lesssim \nu^{-\frac{1}{2}} |\ell|^{-\frac{1}{6}} E_\ell[\omega_\ell]      \cdot  E_{k-\ell}[\theta_{k-\ell}].
		 \end{align}
		 Therefore, combining the above discussion, we obtain a complete estimate for $	R_{1,2,2,2}^M$:
		 \begin{align*}
		 	R_{1,2,2,2}^M \lesssim \nu^{-\frac{1}{2}} \mathcal{E}_{\omega} \cdot \mathcal{E}_{\theta}.
		 \end{align*}
		  Collecting all the estimates for $ R_{1,2}^M$, we then obtain:
		 \begin{align*}
		 		R_{1,2}^M \lesssim \nu^{-\frac{1}{2}} \mathcal{E}_{\omega} \cdot \mathcal{E}_{\theta}.
		 \end{align*}
		 Next, we consider the term $ R_{1,3}^M$. By performing a frequency decomposition, it can be written as three components:
		 \begin{align*}
		 	R_{1,3}^M \lesssim&  \nu^{-\frac{1}{6}}  \left(\int_{0 \le |k-\ell| \le 10 \nu} + \int_{10 \nu \le |k-\ell| \le 1}  + \int_{  |k-\ell| \ge 1}\right) \int_{  |\ell| \ge 1} |k|^{\frac{2}{3}} \left\| u^1_\ell \cdot \theta_{k-\ell}\right\|_{L^2_t L^2_y}   \mathrm{d} \ell  \mathrm{d} k\\
		 	:=& R_{1,3,1}^M   +R_{1,3,2}^M  + R_{1,3,3}^M .
		 \end{align*}
		 For $ R_{1,3,1}^M $ and $ R_{1,3,3}^M $, we can directly apply estimates \eqref{u1_es} and \eqref{theta_es} to obtain:
		 \begin{align*}
		 	R_{1,3,1}^M  &\lesssim \nu^{-\frac{1}{6}}   \int_{0 \le |k-\ell| \le 10 \nu}  \int_{  |\ell| \ge 1} \left(|k-\ell|^{\frac{2}{3}}  +|\ell|^{\frac{2}{3}} \right)\left\| u^1_\ell \right\|_{L^2_t L^\infty_y} \left\| \theta_{k-\ell} \right\|_{L^\infty_t L^2_y}   \mathrm{d} \ell  \mathrm{d} k\\
		 	&\lesssim \nu^{-\frac{7}{24}} \mathcal{E}_{\omega} \cdot \mathcal{E}_{\theta},
		 \end{align*}
		 and
		 \begin{align*}
		 	R_{1,3,3}^M \lesssim& \nu^{-\frac{1}{6}}   \int_{  |k-\ell| \ge 1 }  \int_{  |\ell| \ge 1}  |k-\ell|^{\frac{2}{3}}     \left\| u^1_\ell \right\|_{L^\infty_t L^\infty_y} \left\| \theta_{k-\ell} \right\|_{L^2_t L^2_y}   \mathrm{d} \ell  \mathrm{d} k\\
		 	& + \nu^{-\frac{1}{6}}   \int_{  |k-\ell| \ge 1 }  \int_{  |\ell| \ge 1}  |\ell|^{\frac{2}{3}}  \left\| u^1_\ell \right\|_{L^2_t L^\infty_y} \left\| \theta_{k-\ell} \right\|_{L^\infty_t L^2_y}   \mathrm{d} \ell  \mathrm{d} k \\
\lesssim& \nu^{-\frac{1}{3}} \mathcal{E}_{\omega} \cdot \mathcal{E}_{\theta}.
		 \end{align*}
		 The term $  R_{1,3,2}^M $ requires a further decomposition in order to control the contributions from different frequency interactions separately
		 \begin{align*}
		 R_{1,3,2}^M  \lesssim& \nu^{-\frac{1}{6}}   \int_{10 \nu \le |k-\ell| \le 1}  \int_{  |\ell| \ge 1} \left(|k-\ell|^{\frac{2}{3}}  +|\ell|^{\frac{2}{3}} \right)\left\| u^1_\ell \cdot  \theta_{k-\ell} \right\|_{L^2_t L^2_y}   \mathrm{d} \ell  \mathrm{d} k \\
: =& R_{1,3,2,1}^M  + R_{1,3,2,2}^M .
		 \end{align*}
		 The estimate for the first term can be readily verified as
		 \begin{align*}
		 	R_{1,3,2,1}^M \lesssim  \nu^{-\frac{7}{24}} \mathcal{E}_{\omega} \cdot \mathcal{E}_{\theta}.
		 \end{align*}
		 Analogous to the treatment in estimate \eqref{ta2}, by applying H\"older's inequality with appropriate weights, we obtain
		 \begin{align*}
		 	\nu^{-\frac{1}{6}}   \left\| u^1_\ell \cdot \theta_{k-\ell}\right\|_{L^2_t L^2_y} &\lesssim 	\nu^{-\frac{1}{6}}  \left( \left\| u^1_\ell \right\|_{L^\infty_t L^\infty_y}     \left\|  \theta_{k-\ell} \right\|_{L^2_t L^2_y}  \right)^{\frac{1}{3}}  \cdot \left(  \left\| u^1_\ell \right\|_{L^2_t L^\infty_y}     \left\|  \theta_{k-\ell} \right\|_{L^\infty_t L^2_y}    \right)^{\frac{2}{3}}\\
		 	&\lesssim \nu^{-\frac{5}{12}} |\ell|^{-\frac{2}{3}} E_\ell[\omega_\ell]      \cdot  E_{k-\ell}[\theta_{k-\ell}].
		 \end{align*}
		 From the above, we immediately obtain the estimate for $ R_{1,3,2,2}^M $
		 \begin{align*}
		 	R_{1,3,2,2}^M  \lesssim \nu^{-\frac{5}{12}}  \int_{10 \nu \le |k-\ell| \le 1}  \int_{  |\ell| \ge 1}  |\ell|^{\frac{2}{3}} \cdot  |\ell|^{-\frac{2}{3}} E_\ell[\omega_\ell]      \cdot  E_{k-\ell}[\theta_{k-\ell}] \mathrm{d} \ell  \mathrm{d} k \lesssim \nu^{-\frac{5}{12}} \mathcal{E}_{\omega} \cdot \mathcal{E}_{\theta}.
		 \end{align*}
		 Furthermore, by collecting the estimates for all the decomposed components of $ R_{1,3}^M$, we arrive at the overall estimate for $ R_{1,3}^M$:
		 \begin{align*}
		 	R_{1,3}^M \lesssim \nu^{-\frac{5}{12}} \mathcal{E}_{\omega} \cdot \mathcal{E}_{\theta}.
		 \end{align*}
		 Finally, combining all the estimates for $ R_{1,1}^M$--$ R_{1,3}^M$, we complete the proof of the lemma.
	\end{proof}
\end{lemma}

\subsection{Nonlinear estimates for high frequency}
For the high-frequency regime, specifically, we have
\begin{lemma}\label{prop_theta_high}
	Let $ \theta_k$ be a solution to system \eqref{equation_k}.
	For the high-frequency regime $    |k| \ge 1$, there exists a constant $C$, independent of $\nu$ and $k$, such that
	\begin{align*}
	\left\| |k|^{\frac{1}{3}} \| \theta_k  \|_{L^\infty_t L^2_y} \right\|_{L^1_{|k| \ge 1}} + 	\nu^{\frac{1}{6}} \left\| \| |k|^{\frac{2}{3}} \theta_k  \| _{L^2_t L^2_y} \right\|_{ |k| \ge 1 }  \le C\left\| |k|^{\frac{1}{3}}\theta^{\mathrm{in}}_k \right\|_{L^1_{|k| \ge 1}} + C \nu^{-\frac{1}{2}} \mathcal{E}_{\omega} \cdot \mathcal{E}_{\theta}.
	\end{align*}
	\begin{proof}
		First, we divide the frequencies into two regions, $ \nu^{\frac{1}{6}} |k|^{\frac{1}{3}} \le 1$( corresponding to $ I_4$) and $ \nu^{\frac{1}{6}} |k|^{\frac{1}{3}} \ge 1$( corresponding to $ I_5$), and estimate them separately.
		divide two cases .
	\begin{itemize}
		\item \textbf{Case1:}\ $ \bf \nu^{\frac{1}{6}} |k|^{\frac{1}{3}} \le 1$. In this case, the enhanced dissipation term becomes dominant, thereby exerting primary control over the behavior of the system.
		\\
		
		By a straightforward energy estimate, we obtain
		\begin{align*}
		\big\| E_k[\theta_k ] \big\|_{L^1_{  |k| \ge 1}} &\lesssim \left\| |k|^{\frac{1}{3}} \left\| \theta_k^{\mathrm{in}} \right\|_{L^2_y} \right\|_{L^1_{ |k| \ge 1}} +  \int_{ |k| \ge 1}  \left( \nu^{-\frac{1}{6}} |k|  \left\| g^1_k \right\|_{L^2_t L^2_y} +  \nu^{-\frac{1}{2}}  |k|^{\frac{1}{3}}  \left\| g^2_k \right\|_{L^2_t L^2_y}  \right) \mathrm{d} k\\
		&\lesssim \left\| |k|^{\frac{1}{3}} \left\| \theta_k^{\mathrm{in}} \right\|_{L^2_y} \right\|_{L^1_{|k| \ge 1 }} + R_1^{H,L} + R_2^{H,L}.
		\end{align*}
		Analogous to the treatment in estimate  \eqref{R2_M}, we decompose the relevant term as follows
		\begin{align*}
			R_2^{H,L} 	\lesssim& \nu^{-\frac{1}{2}} \sum_{i = 1}^{3}  \int_{I_3}  \int_{I_i}  |k|^{\frac{1}{3}} \left\|  u^2_\ell   \right\|^{\frac{1}{2}}_{L^2_t L^2_y} \left(|\ell|\left\|  u^1_\ell   \right\|_{L^2_t L^2_y}\right)^{\frac{1}{2}}   \cdot \left\|  \theta_{k-\ell} \right\|_{L^\infty_t L^2_y}  \mathrm{d} \ell  \mathrm{d} k\\
 :=& R_{2,1}^{H,L} + R_{2,2}^{H,L} + R_{2,3}^{H,L}.
		\end{align*}
		For $ 	R_{2,3}^{H,L}$, when $ |k|, |\ell| \ge 1$ and $ |k-\ell| \le 1$, we have $ |k| \approx |\ell|$.
		In this case, by applying estimate \eqref{theta_es}, we obtain
		\begin{align*}
			R_{2,3}^{H,L} &\lesssim \nu^{-\frac{1}{2}}  \left(\int_{0 \le |k-\ell| \le 1}   + \int_{  |k-\ell| \ge 1}\right) \int_{  |\ell| \ge 1} |k|^{\frac{1}{3}} \left\|  u^2_\ell   \right\|^{\frac{1}{2}}_{L^2_t L^2_y} \left(|\ell|\left\|  u^1_\ell   \right\|_{L^2_t L^2_y}\right)^{\frac{1}{2}}   \cdot \left\|  \theta_{k-\ell} \right\|_{L^\infty_t L^2_y}   \mathrm{d} \ell  \mathrm{d} k\\
			&\lesssim \nu^{-\frac{1}{2}} \int_{0 \le |k-\ell| \le 1} \int_{  |\ell| \ge 1} |k|^{\frac{1}{3}} |\ell|^{-\frac{1}{2}}   E_\ell[\omega_\ell]       \cdot  E_{k-\ell}[\theta_{k-\ell}]   \mathrm{d} \ell  \mathrm{d} k\\
			&\quad + \nu^{-\frac{1}{2}}  \int_{  |k-\ell| \ge 1} \int_{  |\ell| \ge 1} \left(|k-\ell|^{\frac{1}{3}} + |\ell|^{\frac{1}{3}}\right) |\ell|^{-\frac{1}{2}}   E_\ell[\omega_\ell]      \cdot \left\|  \theta_{k-\ell} \right\|_{L^\infty_t L^2_y}   \mathrm{d} \ell  \mathrm{d} k\\
			&\lesssim \nu^{-\frac{1}{2}} \mathcal{E}_{\omega} \cdot \mathcal{E}_{\theta}.
		\end{align*}
		For $	R_{2,1}^{H,L}$, when  $ |k| \ge 1$  and $ |\ell| \le 10 \nu$, since $10\nu \le 10\nu_0 \le \frac{1}{2}$, we have $|k-\ell| \approx |k|$.
		Combining estimate \eqref{theta_es} with the Poincar\'e inequality, we obtain
		\begin{align*}
			R_{2,1}^{H,L}  &\lesssim \nu^{-\frac{1}{2}}   \int_{|k| \ge 1}  \int_{|\ell| \le 10 \nu}  |k|^{\frac{1}{3}} \left\|  u^2_\ell   \right\|^{\frac{1}{2}}_{L^2_t L^2_y} \left(|\ell|\left\|  u^1_\ell   \right\|_{L^2_t L^2_y}\right)^{\frac{1}{2}}   \cdot \left\|  \theta_{k-\ell} \right\|_{L^\infty_t L^2_y}  \mathrm{d} \ell  \mathrm{d} k \\
			&\lesssim \nu^{-\frac{1}{2}}   \int_{|k-\ell| \gtrsim 1}  \int_{|\ell| \le 10 \nu}  |k-\ell|^{\frac{1}{3}}  |\ell|\left\|  u^1_\ell   \right\|_{L^2_t L^2_y}    \cdot \left\|  \theta_{k-\ell} \right\|_{L^\infty_t L^2_y}  \mathrm{d} \ell  \mathrm{d} k \\
			&\lesssim \nu^{-\frac{1}{2}} \mathcal{E}_{\omega} \cdot \mathcal{E}_{\theta}.
		\end{align*}
		For $R_{2,2}^{H,L}$, when $ |k| \ge 1$,  $  10 \nu \le |\ell| \le 1$ and $ |k-\ell| \ge 1$, we have $|k-\ell| \approx |k|$.
		By an argument analogous to that for $	R_{2,1}^{H,L} $, we obtain
		\begin{align*}
			 R_{2,2}^{H,L}    \lesssim& \nu^{-\frac{1}{2}}  \int_{0 \le |k-\ell| \le 1}   \int_{ 10 \nu \le  |\ell| \le 1} |k|^{\frac{1}{3}} \left\|  u^2_\ell   \right\|^{\frac{1}{2}}_{L^2_t L^2_y} \left(|\ell|\left\|  u^1_\ell   \right\|_{L^2_t L^2_y}\right)^{\frac{1}{2}}   \cdot \left\|  \theta_{k-\ell} \right\|_{L^\infty_t L^2_y}   \mathrm{d} \ell  \mathrm{d} k\\
&+\nu^{-\frac{1}{2}}  \int_{  |k-\ell| \ge 1} \int_{ 10 \nu \le  |\ell| \le 1} |k|^{\frac{1}{3}} \left\|  u^2_\ell   \right\|^{\frac{1}{2}}_{L^2_t L^2_y} \left(|\ell|\left\|  u^1_\ell   \right\|_{L^2_t L^2_y}\right)^{\frac{1}{2}}   \cdot \left\|  \theta_{k-\ell} \right\|_{L^\infty_t L^2_y}   \mathrm{d} \ell  \mathrm{d} k\\
			 \lesssim&   \nu^{-\frac{1}{2}} \int_{0 \le |k-\ell| \le 1} \int_{  10 \nu \le |\ell| \le 1} \left( |k-\ell|^{\frac{1}{3}} + |\ell|^{\frac{1}{3}}\right)   E_\ell[\omega_\ell]       \cdot  E_{k-\ell}[\theta_{k-\ell}]   \mathrm{d} \ell  \mathrm{d} k\\
			 &+  \nu^{-\frac{1}{2}}  \int_{  |k-\ell| \ge 1} \int_{ 10 \nu \le  |\ell| \le 1}  |k-\ell|^{\frac{1}{3}}     E_\ell[\omega_\ell]      \cdot \left\|  \theta_{k-\ell} \right\|_{L^\infty_t L^2_y}   \mathrm{d} \ell  \mathrm{d} k \\
\lesssim& \nu^{-\frac{1}{2}} \mathcal{E}_{\omega} \cdot \mathcal{E}_{\theta}.
		\end{align*}
		By combining the estimates for $R_{2,1}^{H,L} $ through $R_{2,3}^{H,L} $, we obtain
		\begin{align*}
			R_{2}^{H,L}  \lesssim \nu^{-\frac{1}{2}} \mathcal{E}_{\omega} \cdot \mathcal{E}_{\theta}.
		\end{align*}
		Next, we turn to the most delicate term, $	R_1^{H,L}$.
		To handle it, we first further decompose it into the following three parts:
			\begin{align*}
				R_1^{H,L}  \lesssim \nu^{-\frac{1}{6}} \sum_{i=1}^{3} \int_{ |k| \ge 1} \int_{I_i}    |k|  \left\| u^1_\ell \cdot \theta_{k-\ell} \right\|_{L^2_t L^2_y}  \mathrm{d} \ell \mathrm{d} k := R_{1,1}^{H,L} + R_{1,2}^{H,L} + R_{1,3}^{H,L}.
			\end{align*}
			Each of these components will be estimated separately according to its respective frequency structure.
			For $R_{1,1}^{H,L}$, when $|k| \ge 1$ and $|\ell| \le 10\nu$, the condition $10\nu \le 10\nu_{0} \le \tfrac{1}{2}$ implies that $|k-\ell| \approx |k|$.
			Using this observation together with the corresponding estimates in \eqref{u1_es} and \eqref{theta_es}, we can further derive
			\begin{align*}
				R_{1,1}^{H,L}  &\lesssim \nu^{-\frac{1}{6}} \int_{ |k-\ell| \gtrsim 1} \int_{ |\ell| \le 10 \nu}  |k|^{\frac{1}{3}}  \left\| u^1_\ell \right\|_{L^\infty_t L^\infty_y}      |k-\ell|^{\frac{2}{3}}  \left\|  \theta_{k-\ell} \right\|_{L^2_t L^2_y}   \mathrm{d} \ell \mathrm{d} k \\
				&\lesssim \nu^{-\frac{1}{2}} \int_{ |k-\ell| \gtrsim 1} \int_{ |\ell| \le 10 \nu}  \nu^{\frac{1}{6}}|k|^{\frac{1}{3}}      E_\ell[\omega_\ell]       \cdot  E_{k-\ell}[\theta_{k-\ell}]       \mathrm{d} \ell \mathrm{d} k\\
				&\lesssim  \nu^{-\frac{1}{2}} \mathcal{E}_{\omega} \cdot \mathcal{E}_{\theta},
			\end{align*}
			where, in the last step of the above estimate, we have used the fact that $\nu^{\frac{1}{6}} |k|^{\frac{1}{3}} \le 1$.

			For $	R_{1,2}^{H,L}$, when the high-frequency modes ($|k| \ge 1$) interact with the medium-frequency modes ($10\nu \le |\ell| \le 1$), the frequency ranges still exhibit a certain degree of complexity.
			Therefore, we further decompose $	R_{1,2}^{H,L}$ into three parts, which will be estimated separately:
			\begin{align}\label{ta3}
					\qquad R_{1,2}^{H,L} \lesssim& \nu^{-\frac{1}{6}}\left( \int_{0 \le |k-\ell| \le 10 \nu} + \int_{10 \nu \le |k-\ell| \le 1} +  \int_{ |k-\ell| \ge 1} \right)\int_{  10 \nu \le |\ell| \le 1}  |k|  \left\| u^1_\ell \cdot \theta_{k-\ell} \right\|_{L^2_t L^2_y}     \mathrm{d} \ell \mathrm{d} k \nn\\
					 : =& 	R_{1,2,1}^{H,L} + 	R_{1,2,2}^{H,L} + 	R_{123}^{H,L}.
			\end{align}
			For the first term in \eqref{ta3}, within the current frequency range we observe that
			\[
			\frac{|k|}{|\ell|^{3/8}}
			\lesssim |\ell|^{5/8} + \frac{|k-\ell|}{|\ell|^{3/8}}
			\lesssim 1 .
			\]
			Hence, this term is well controlled in this regime.
			Combining this with the estimates in \eqref{u1_es} and \eqref{theta_es}, we obtain
			\begin{align*}
				R_{1,2,1}^{H,L} \lesssim \nu^{-\frac{1}{6}} \int_{0 \le |k-\ell| \le 10 \nu} \int_{ 10 \nu \le |\ell| \le 1}  |k|  \left\| u^1_\ell \right\|_{L^2_t L^\infty_y} \left\| \theta_{k-\ell} \right\|_{L^\infty_t L^2_y}       \mathrm{d} \ell \mathrm{d} k \lesssim \nu^{-\frac{7}{24}} \mathcal{E}_{\omega} \cdot \mathcal{E}_{\theta}.
			\end{align*}
			For the second term in \eqref{ta3}, when $|k| \ge 1$, $10\nu \le |\ell| \le 1$, and $10\nu \le |k-\ell| \le 1$,
			we observe that $k$ and $\ell$ must share the same sign.
			Consequently, by applying the triangle inequality, we further decompose this term into the following two parts:
			\begin{align*}
				R_{1,2,2}^{H,L} \lesssim& \nu^{-\frac{1}{6}} \int_{ 10 \nu \le |k-\ell| \le 1} \int_{ 10 \nu \le |\ell| \le 1}  \left(|k-\ell| + |\ell| \right)  \left\| u^1_\ell \cdot \theta_{k-\ell} \right\|_{L^2_t L^2_y}       \mathrm{d} \ell \mathrm{d} k\\
 : =& R_{1,2,2,1}^{H,L} + R_{1,2,2,2}^{H,L}.
			\end{align*}
			For the term $R_{1,2,2,1}^{H,L} $, within the current frequency range we observe that
			\[
			|\ell|^{-3/16} \lesssim \nu^{-3/16}.
			\]
			Therefore, by applying H\"older's inequality with suitable weights, we obtain
			\begin{align*}
				R_{1,2,2,1}^{H,L} &\lesssim \nu^{-\frac{1}{6}} \int_{ 10 \nu \le |k-\ell| \le 1} \int_{ 10 \nu \le |\ell| \le 1}   |k-\ell|   \left( \left\| u^1_\ell \right\|_{L^\infty_t L^\infty_y}     \left\|  \theta_{k-\ell} \right\|_{L^2_t L^2_y}  \right)^{\frac{1}{2}} \\
 &\qquad\qquad\qquad\qquad\qquad\qquad\qquad\quad\times \left(  \left\| u^1_\ell \right\|_{L^2_t L^\infty_y}     \left\|  \theta_{k-\ell} \right\|_{L^\infty_t L^2_y}    \right)^{\frac{1}{2}}     \mathrm{d} \ell \mathrm{d} k\\
				  &\lesssim \nu^{-\frac{1}{2}} \int_{ 10 \nu \le |k-\ell| \le 1} \int_{ 10 \nu \le |\ell| \le 1}  |k-\ell| |k-\ell|^{-\frac{1}{6}}   E_\ell[\omega_\ell]      \cdot  E_{k-\ell}[\theta_{k-\ell}]        \mathrm{d} \ell \mathrm{d} k \\
				   &\lesssim \nu^{-\frac{1}{2 }} \mathcal{E}_{\omega} \cdot \mathcal{E}_{\theta}.
			\end{align*}
			When the frequencies satisfy $10\nu \le |k-\ell| \le 1$ and $10\nu \le |\ell| \le 1$, we have
			\[
			|\ell|^{-5/27} \lesssim \nu^{-5/27}.
			\]
			Combining this with the estimates in \eqref{u1_es} and \eqref{theta_es}, we obtain
			\begin{align*}
			\nu^{-\frac{1}{6}} \left\| u^1_\ell \cdot \theta_{k-\ell} \right\|_{L^2_t L^2_y}   &\lesssim 	\nu^{-\frac{1}{6}}  \left( \left\| u^1_\ell \right\|_{L^\infty_t L^\infty_y}     \left\|  \theta_{k-\ell} \right\|_{L^2_t L^2_y}  \right)^{\frac{5}{9}}  \cdot \left(  \left\| u^1_\ell \right\|_{L^2_t L^\infty_y}     \left\|  \theta_{k-\ell} \right\|_{L^\infty_t L^2_y}    \right)^{\frac{4}{9}}\\
			&\lesssim \nu^{-\frac{1}{ 2}} |\ell|^{-\frac{1}{6}} E_\ell[\omega_\ell]      \cdot  E_{k-\ell}[\theta_{k-\ell}].
			\end{align*}
			Based on the above estimates, we conclude that $	R_{1,2,2,2}^{H,L}$ satisfies the following bound:
			\begin{align*}
			R_{1,2,2,2}^{H,L}  \lesssim \nu^{-\frac{1}{2}} \int_{ 10 \nu \le |k-\ell| \le 1} \int_{ 10 \nu \le |\ell| \le 1}   |\ell|^{\frac{5}{6}}   E_\ell[\omega_\ell]      \cdot  E_{k-\ell}[\theta_{k-\ell}]        \mathrm{d} \ell \mathrm{d} k  \lesssim \nu^{-\frac{1}{2 }} \mathcal{E}_{\omega} \cdot \mathcal{E}_{\theta}.
			\end{align*}
			For the third term in \eqref{ta3}, when $|k| \ge 1$, $10\nu \le |\ell| \le 1$, and $|k-\ell| \ge 1$,
			since $|k-\ell| \approx |k|$, we can directly apply the estimates in \eqref{u1_es} and \eqref{theta_es} to obtain
			\begin{align*}
				R_{123}^{H,L} \lesssim& \nu^{-\frac{1}{6}} \int_{  |k-\ell| \ge 1} \int_{ 10 \nu \le |\ell| \le 1}  |k|^{\frac{1}{3}}  \left\| u^1_\ell \right\|_{L^\infty_t L^\infty_y} \cdot |k-\ell|^{\frac{2}{3}}  \left\| \theta_{k-\ell} \right\|_{L^2_t L^2_y}       \mathrm{d} \ell \mathrm{d} k \\
 \lesssim& \nu^{-\frac{1}{2 }} \mathcal{E}_{\omega} \cdot \mathcal{E}_{\theta},
			\end{align*}
			where	in the final step of the above estimate, we have used the fact that $\nu^{\frac{1}{6}} |k|^{\frac{1}{3}} \le 1$.

			By combining the estimates for the various components of $	R_{1,2}^{H,L} $ discussed above, we obtain
			\begin{align*}
				R_{1,2}^{H,L}  \lesssim \nu^{-\frac{1}{2 }} \mathcal{E}_{\omega} \cdot \mathcal{E}_{\theta}.
			\end{align*}
			For the estimate of $	R_{1,3}^{H,L}$, we similarly employ a frequency decomposition method analogous to that in \eqref{ta3}, yielding
			\begin{align*}
			R_{1,3}^{H,L}  \lesssim& \nu^{-\frac{1}{6}}\left( \int_{0 \le |k-\ell| \le 10 \nu} + \int_{10 \nu \le |k-\ell| \le 1} +  \int_{ |k-\ell| \ge 1} \right)\int_{    |\ell| \ge 1}  |k|  \left\| u^1_\ell \cdot \theta_{k-\ell} \right\|_{L^2_t L^2_y}     \mathrm{d} \ell \mathrm{d} k \\
			 : =& 	R_{1,3,1}^{H,L} + 	R_{1,3,2}^{H,L} + 	R_{1,3,3}^{H,L}.
			\end{align*}
			This allows us to uniformly control the contributions of $	R_{1,3}^{H,L}$ across the different frequency regimes. For $R_{1,3,1}^{H,L}$ and $R_{1,3,2}^{H,L}$, when $0 \le |k-\ell| \le 1$, $|k| \ge 1$, and $|\ell| \ge 1$, we have$	|k| \approx |\ell|$
			within this frequency range.
			Therefore, we can directly apply \eqref{u1_es} and \eqref{theta_es} to obtain the corresponding estimates:
			\begin{align*}
				R_{1,3,1}^{H,L} + R_{1,3,2}^{H,L} &\lesssim \nu^{-\frac{1}{6}} \int_{0 \le |k-\ell| \le 1} \int_{    |\ell| \ge 1}  |k|   \left\| u^1_\ell \right\|_{L^2_t L^\infty_y}     \left\|  \theta_{k-\ell} \right\|_{L^\infty_t L^2_y}     \mathrm{d} \ell \mathrm{d} k\\
				&\lesssim \nu^{-\frac{7}{24}} \int_{0 \le |k-\ell| \le 1} \int_{    |\ell| \ge 1}  |k|^{\frac{1}{4}} E_\ell[\omega_\ell]      \cdot  E_{k-\ell}[\theta_{k-\ell}]    \mathrm{d} \ell \mathrm{d} k\\
				&\lesssim \nu^{-\frac{5}{12}} \mathcal{E}_{\omega} \cdot \mathcal{E}_{\theta}.
			\end{align*}
			For the term $R_{1,3,3}^{H,L}$, we first apply Lemma \ref{lebasis2} to decompose it into the following two parts:
			\begin{align*}
			R_{1,3,3}^{H,L} \lesssim& \nu^{-\frac{1}{6}} \int_{  |k-\ell| \ge 1} \int_{    |\ell| \ge 1}  \left( |k|^m \cdot |k-\ell|^{1-m} + |k|^n \cdot |\ell|^{1-n}\right)   \left\| u^1_\ell  \cdot \theta_{k-\ell} \right\|_{L^2_t L^2_y}     \mathrm{d} \ell \mathrm{d} k\\
			:=& R_{1,3,3,1}^{H,L} + R_{1,3,3,2}^{H,L}.
			\end{align*}
			For the first part, $	R_{1,3,3,1}^{H,L}$, taking $m = \frac{1}{3}$ and using the key fact $\nu^{\frac{1}{6}} |k|^{\frac{1}{3}} \le 1$, together with the estimates in \eqref{u1_es} and \eqref{theta_es}, we obtain
			\begin{align*}
			R_{1,3,3,1}^{H,L} \lesssim& \nu^{-\frac{1}{6}} \int_{  |k-\ell| \ge 1} \int_{    |\ell| \ge 1}   |k|^{\frac{1}{3}} \cdot |k-\ell|^{\frac{2}{3}} \left\| u^1_\ell \right\|_{L^\infty_t L^\infty_y}     \left\|  \theta_{k-\ell} \right\|_{L^2_t L^2_y} \mathrm{d} \ell \mathrm{d} k\\
			\lesssim& \nu^{-\frac{1}{3}} \int_{  |k-\ell| \ge 1} \int_{    |\ell| \ge 1}   |k|^{\frac{1}{3}} |\ell|^{-\frac{1}{2}} E_\ell[\omega_\ell]      \cdot  E_{k-\ell}[\theta_{k-\ell}]    \mathrm{d} \ell \mathrm{d} k\\
 \lesssim& \nu^{-\frac{1}{ 2}} \mathcal{E}_{\omega} \cdot \mathcal{E}_{\theta}.
			\end{align*}
			For the second part, $R_{1,3,3,2}^{H,L}$, in the frequency range where $|k-\ell| \ge 1$ and $|\ell| \ge 1$, by applying the weighted H\"older's inequality, we have
			\begin{align*}
			\nu^{-\frac{1}{6}} \left\| u^1_\ell \cdot \theta_{k-\ell} \right\|_{L^2_t L^2_y}   &\lesssim 	\nu^{-\frac{1}{6}}  \left( \left\| u^1_\ell \right\|_{L^\infty_t L^\infty_y}     \left\|  \theta_{k-\ell} \right\|_{L^2_t L^2_y}  \right)^{\frac{1}{3}}  \cdot \left(  \left\| u^1_\ell \right\|_{L^2_t L^\infty_y}     \left\|  \theta_{k-\ell} \right\|_{L^\infty_t L^2_y}    \right)^{\frac{2}{3}}\\
			&\lesssim \nu^{-\frac{11}{ 36}} |k-\ell|^{-\frac{1}{3}} |\ell|^{-\frac{2}{3}} E_\ell[\omega_\ell]      \cdot  E_{k-\ell}[\theta_{k-\ell}] .
			\end{align*}
			Furthermore, taking $n = \frac{1}{3}$ and again using $\nu^{\frac{1}{6}} |k|^{\frac{1}{3}} \le 1$, we can deduce
			\begin{align*}
			R_{1,3,3,2}^{H,L} \lesssim& \nu^{-\frac{11}{36}} \int_{  |k-\ell| \ge 1} \int_{    |\ell| \ge 1}   |k|^{\frac{1}{3}} \cdot |\ell|^{\frac{2}{3}}  |k-\ell|^{-\frac{1}{3}} |\ell|^{-\frac{2}{3}} E_\ell[\omega_\ell]      \cdot  E_{k-\ell}[\theta_{k-\ell}]  \mathrm{d} \ell \mathrm{d} k\\
			\lesssim& \nu^{-\frac{17}{36}} \int_{  |k-\ell| \ge 1} \int_{    |\ell| \ge 1}     E_\ell[\omega_\ell]      \cdot  E_{k-\ell}[\theta_{k-\ell}]    \mathrm{d} \ell \mathrm{d} k\\
 \lesssim& \nu^{-\frac{17}{36}} \mathcal{E}_{\omega} \cdot \mathcal{E}_{\theta}.
			\end{align*}
			By combining all the above estimates for the components of $R_{1,3}^{H,L}$, we arrive at
			\begin{align*}
			R_{1,3}^{H,L} \lesssim  \nu^{-\frac{1}{ 2}} \mathcal{E}_{\omega} \cdot \mathcal{E}_{\theta}.
			\end{align*}
			By combining the estimates for $ R_{1,1}^{H,L}$ through $ R_{1,3}^{H,L}$, one can get
			\begin{align*}
			R_{1}^{H,L} \lesssim  \nu^{-\frac{1}{ 2}} \mathcal{E}_{\omega} \cdot \mathcal{E}_{\theta}.
			\end{align*}
			 We thus complete the proof of the lemma in the case $\nu^{\frac{1}{6}} |k|^{\frac{1}{3}} \le 1$.

		\item \textbf{Case2:}\ $ \bf \nu^{\frac{1}{6}} |k|^{\frac{1}{3}} \ge 1$. In this situation, Lemma \ref{le_theta} is no longer sufficient to provide further information on  dissipation, and the viscous term becomes the dominant factor. \\
		
		Using the definition of $E_k$ and standard energy estimates, together with Lemma  \ref{lebasis3}, we obtain
			\begin{align*}
			\big\| E_k[\theta_k ] \big\|_{L^1_{  |k| \ge 1}} &\lesssim \left\| |k|^{\frac{1}{3}} \left\| \theta_k^{\mathrm{in}} \right\|_{L^2_y} \right\|_{L^1_{ |k| \ge 1}} +  \int_{ |k| \ge 1}  \nu^{-\frac{1}{2}} |k|^{\frac{1}{3}} \left(   \left\| g^1_k \right\|_{L^2_t L^2_y} +     \left\| g^2_k \right\|_{L^2_t L^2_y}  \right) \mathrm{d} k\\
			&\lesssim \left\| |k|^{\frac{1}{3}} \left\| \theta_k^{\mathrm{in}} \right\|_{L^2_y} \right\|_{L^1_{|k| \ge 1 }} + R_1^{H,H} + R_2^{H,H}.
			\end{align*}
			For the term $ R_2^{H,H} $, we can proceed in a manner similar to that for $ R_2^{H,L}$, yielding
			\begin{align*}
			R_2^{H,H} \lesssim \nu^{-\frac{1}{ 2}} \mathcal{E}_{\omega} \cdot \mathcal{E}_{\theta}.
			\end{align*}
			Therefore, we omit further details for this estimate.
				Next, we focus on the term $R_1^{H,H} $. Using a decomposition analogous to that in \eqref{ta3}, we have
			\begin{align*}
			R_1^{H,H}  \lesssim \nu^{-\frac{1}{2}} \sum_{i=1}^{3} \int_{ |k| \ge 1} \int_{I_i}    |k|^{\frac{1}{3}}  \left\| u^1_\ell \cdot \theta_{k-\ell} \right\|_{L^2_t L^2_y}  \mathrm{d} \ell \mathrm{d} k := R_{1,1}^{H,H} + R_{1,2}^{H,H} + R_{1,3}^{H,H}.
			\end{align*}
			For the term $	R_{1,1}^{H,H}$, when $|k|\ge 1$ and $|\ell|\le 10\nu$, we have the frequency comparability $ |k-\ell| \approx |k|$
			Hence, combining the estimates \eqref{u1_es} and \eqref{theta_es}, we obtain
			\begin{align*}
				R_{1,1}^{H,H} \lesssim&\nu^{-\frac{1}{2}}   \int_{ |k-\ell| \gtrsim 1} \int_{I_1}    |k|^{\frac{1}{3}} \left\| u^1_\ell \right\|_{L^\infty_t L^\infty_y}     \left\|  \theta_{k-\ell} \right\|_{L^2_t L^2_y}    \mathrm{d} \ell \mathrm{d} k \\
				\lesssim& \nu^{-\frac{2}{3}}   \int_{ |k-\ell| \gtrsim 1} \int_{I_1}    |k|^{\frac{1}{3}} \cdot   |k-\ell|^{-\frac{2}{3}} E_\ell[\omega_\ell]      \cdot  E_{k-\ell}[\theta_{k-\ell}]     \mathrm{d} \ell \mathrm{d} k \\
\lesssim& \nu^{-\frac{1}{ 2}} \mathcal{E}_{\omega} \cdot \mathcal{E}_{\theta},
			\end{align*}
		where have we used the fact $ \nu^{-\frac{1}{6}} |k|^{-\frac{1}{3}} \le 1 .$
			For the term $	R_{1,2}^{H,H}  $, since it involves the interaction between high and intermediate frequencies, its structure is more delicate.
			Following the same decomposition strategy as in \eqref{ta3}, we further split $	R_{1,2}^{H,H}  $ into the following three parts:
			\begin{align*}
				\qquad R_{1,2}^{H,H}    \lesssim& \nu^{-\frac{1}{2}}\left( \int_{0 \le |k-\ell| \le 10 \nu} + \int_{10 \nu \le |k-\ell| \le 1} +  \int_{ |k-\ell| \ge 1} \right)\int_{    10 \nu \le |\ell| \le 1}  |k|^{\frac{1}{3}}  \left\| u^1_\ell \cdot \theta_{k-\ell} \right\|_{L^2_t L^2_y}     \mathrm{d} \ell \mathrm{d} k \\
				: =& 	R_{1,2,1}^{H,H} + 	R_{1,2,2}^{H,H} + 	R_{123}^{H,H}.
			\end{align*}
			For the term $	R_{1,2,1}^{H,H}$, when $|k| \ge 1$, $10\nu \le |\ell| \le 1$ and $|k-\ell| \le 10\nu$, the condition
			$10\nu \le 10\nu_{0} \le \tfrac{1}{2}$ implies that
			\[
			|k| \approx |\ell|.
			\]
			Therefore, applying the estimates in \eqref{u1_es} and \eqref{theta_es}, we obtain
			\begin{align*}
			R_{1,2,1}^{H,H}  \lesssim& \nu^{-\frac{1}{2}}   \int_{ |k-\ell| \le 10 \nu} \int_{|\ell| \gtrsim 1}    |k|^{\frac{1}{3}} \left\| u^1_\ell \right\|_{L^2_t L^\infty_y}     \left\|  \theta_{k-\ell} \right\|_{L^\infty_t L^2_y}    \mathrm{d} \ell \mathrm{d} k \\
			\lesssim& \nu^{-\frac{1}{2}}   \int_{ |k-\ell| \le 10 \nu} \int_{|\ell| \gtrsim 1}    |k|^{\frac{1}{3}} \nu^{-\frac{1}{8}}  |\ell|^{-\frac{3}{4}} E_\ell[\omega_\ell]      \cdot  E_{k-\ell}[\theta_{k-\ell}]     \mathrm{d} \ell \mathrm{d} k \\
\lesssim& \nu^{-\frac{1}{ 2}} \mathcal{E}_{\omega} \cdot \mathcal{E}_{\theta},
			\end{align*}
			where in the final step we have used the fact that  $ \nu^{-\frac{1}{8}}  |k|^{\frac{1}{3}}  |\ell|^{-\frac{3}{4}} \lesssim \nu^{-\frac{1}{8}} |k|^{-\frac{5}{12}} \lesssim 1.$

			For the term $R_{1,2,2}^{H,H}$, we further decompose the frequency of $l$ into the following two parts:
			\begin{align*}
			R_{1,2,2}^{H,H} =& \nu^{-\frac{1}{2}}  \int_{10 \nu \le |k-\ell| \le 1} \left(\int_{    10 \nu \le |\ell| \le \frac{1}{2}}  + \int_{    \frac{1}{2} \le |\ell| \le 1}  \right) |k|^{\frac{1}{3}}  \left\| u^1_\ell \cdot \theta_{k-\ell} \right\|_{L^2_t L^2_y}     \mathrm{d} \ell \mathrm{d} k \\
: =& R_{1,2,2,1}^{H,H} + R_{1,2,2,2}^{H,H}.
			\end{align*}
			For the second part in the above decomposition, the frequency conditions imply that
			\[
			|k| \ge 1,\qquad \tfrac{1}{2} \le |\ell| \le 1,\qquad 10\nu \le |k-\ell| \le 1,
			\]
			and therefore
			$
			|k| \approx |\ell|.
			$
			Proceeding as in the estimate of $R_{1,2,1}^{H,H}$, and using \eqref{u1_es}-\eqref{theta_es}, we obtain
			$ R_{1,2,2,2}^{H,H} \lesssim \nu^{-\frac{1}{ 2}} \mathcal{E}_{\omega} \cdot \mathcal{E}_{\theta}.$
			
			For the first part in the decomposition of $R_{1,2,2}^{H,H} $, when
			$$
			|k| \ge 1,\qquad 10\nu \le |\ell| \le \tfrac{1}{2},\qquad 10\nu \le |k-\ell| \le 1,
			$$
			the condition $10\nu \le 10\nu_0 \le \tfrac12$ yields
		$
			|k| \approx |k-\ell|.
			$
			Hence, we can directly apply the estimates \eqref{u1_es}-\eqref{theta_es} to deduce that
			\begin{align*}
				R_{1,2,2,1}^{H,H} \lesssim& \nu^{-\frac{1}{2}} \int_{  |k-\ell| \gtrsim 1}  \int_{    10 \nu \le |\ell| \le \frac{1}{2}} |k-\ell|^{\frac{1}{3}} \cdot  \nu^{-\frac{1}{6}} |k-\ell|^{-\frac{2}{3}}    E_\ell[\omega_\ell]      \cdot  E_{k-\ell}[\theta_{k-\ell}]        \mathrm{d} \ell \mathrm{d} k\\
 \lesssim& \nu^{-\frac{1}{ 2}} \mathcal{E}_{\omega} \cdot \mathcal{E}_{\theta},
			\end{align*}
			where we have used the fact $ \nu^{-\frac{1}{6}} |k-\ell|^{-\frac{1}{3}} \lesssim 1.$\\
			For the term $R_{123}^{H,H}$, we observe that under the frequency conditions
			\[
			|k| \ge 1,\qquad 10\nu \le |\ell| \le \tfrac12,\qquad |k-\ell| \ge 1,
			\]
			we have
			$
			|k| \approx |k-\ell|.
			$
			Hence, this term can be treated in the same manner as $R_{1,2,2,1}^{H,H}$, which yields
			\begin{align*}
			R_{123}^{H,H} \lesssim \nu^{-\frac{1}{ 2}} \mathcal{E}_{\omega} \cdot \mathcal{E}_{\theta}.
			\end{align*}
			Combining the estimates for all components in $R_{1,2}^{H,H}$, we obtain
			\begin{align*}
			R_{1,2}^{H,H} \lesssim \nu^{-\frac{1}{ 2}} \mathcal{E}_{\omega} \cdot \mathcal{E}_{\theta}.
			\end{align*}
			Next, we turn to the estimate of $R_{1,3}^{H,H}$. Following the same frequency decomposition as in \eqref{ta3}, we split it into the interactions among low, medium, and high frequencies with the high-frequency part:
			\begin{align*}
			R_{1,3}^{H,H}   \lesssim& \nu^{-\frac{1}{2}}\left( \int_{0 \le |k-\ell| \le 10 \nu} + \int_{10 \nu \le |k-\ell| \le 1} +  \int_{ |k-\ell| \ge 1} \right)\int_{     |\ell| \ge 1}  |k|^{\frac{1}{3}}  \left\| u^1_\ell \cdot \theta_{k-\ell} \right\|_{L^2_t L^2_y}     \mathrm{d} \ell \mathrm{d} k \\
			: =& 	R_{1,3,1}^{H,H} + 	R_{1,3,2}^{H,H} + 	R_{1,3,3}^{H,H}.
			\end{align*}
			For the terms $R_{1,3,1}^{H,H}$ and $R_{1,3,2}^{H,H}$, under the frequency regime
			\[
			|k| \ge 1,\qquad |\ell| \ge 1,\qquad 0 \le |k-\ell| \le 1,
			\]
			we have
			$
			|k| \approx |\ell|.
			$
			Therefore, by an argument identical to that used for $R_{1,2,1}^{H,H}$, we directly obtain
			\begin{align*}
			R_{1,3,1}^{H,H} + 	R_{1,3,2}^{H,H} \lesssim \nu^{-\frac{1}{ 2}} \mathcal{E}_{\omega} \cdot \mathcal{E}_{\theta}.
			\end{align*}
			
			For the term $R_{1,3,3}^{H,H}$, we apply Lemma \ref{lebasis2} to transfer the frequency $k$ to the frequencies $\ell$ and $k-\ell$, thereby decomposing it into the following two parts:
			\begin{align*}
			R_{1,3,3}^{H,H} \lesssim& \nu^{-\frac{1}{2}} \int_{  |k-\ell| \ge 1} \int_{    |\ell| \ge 1}  \left( |k|^m \cdot |k-\ell|^{\frac{1}{3}-m} + |k|^n \cdot |\ell|^{\frac{1}{3}-n}\right)   \left\| u^1_\ell  \cdot \theta_{k-\ell} \right\|_{L^2_t L^2_y}     \mathrm{d} \ell \mathrm{d} k\\
			:=& R_{1331}^{H,H} + R_{1332}^{H,H}
			\end{align*}
			For the first part, $R_{1331}^{H,H} $, taking $m = \frac{1}{6}$ and using the estimates in \eqref{u1_es} and \eqref{theta_es}, we obtain
			\begin{align*}
			R_{1331}^{H,H} \lesssim& \nu^{-\frac{1}{2}} \int_{  |k-\ell| \ge 1} \int_{    |\ell| \ge 1}   |k|^{\frac{1}{6}} \cdot |k-\ell|^{\frac{1}{6}} \left\| u^1_\ell \right\|_{L^\infty_t L^\infty_y}     \left\|  \theta_{k-\ell} \right\|_{L^2_t L^2_y} \mathrm{d} \ell \mathrm{d} k\\
			\lesssim& \nu^{-\frac{1}{2}} \int_{  |k-\ell| \ge 1} \int_{    |\ell| \ge 1}   |k|^{\frac{1}{6}} |k-\ell|^{\frac{1}{6}} |\ell|^{-\frac{1}{2}} \cdot \nu^{-\frac{1}{6}} |k-\ell|^{-\frac{2}{3}} E_\ell[\omega_\ell]      \cdot  E_{k-\ell}[\theta_{k-\ell}]    \mathrm{d} \ell \mathrm{d} k \\
\lesssim& \nu^{-\frac{1}{ 2}} \mathcal{E}_{\omega} \cdot \mathcal{E}_{\theta}.
			\end{align*}
			In the final step, we use the key fact:
			$$\nu^{-\frac{1}{6}} |\ell|^{-\frac{1}{2}} |k-\ell|^{-\frac{1}{2}}   |k|^{\frac{1}{6}} \lesssim \nu^{-\frac{1}{6}} |\ell|^{-\frac{1}{2}} |k-\ell|^{-\frac{1}{2}}   |k-\ell|^{\frac{1}{6}}  |\ell|^{\frac{1}{6}} \lesssim   \nu^{-\frac{1}{6}}|k|^{-\frac{1}{3}}  \le 1. $$
			For the second part, $R_{1332}^{H,H}$, taking $n = \frac{1}{12}$, and under the frequency conditions $|k-\ell| \ge 1$ and $|\ell| \ge 1$, we have
			\[
			|\ell|^{-1} |k-\ell|^{-1} \lesssim |k|^{-1}.
			\]
			Hence, by combining the estimates in \eqref{u1_es}-\eqref{theta_es}, we obtain the control
			\begin{align*}
			R_{1332}^{H,H} \lesssim& \nu^{-\frac{1}{2}} \int_{  |k-\ell| \ge 1} \int_{    |\ell| \ge 1}   |k|^{\frac{1}{12}} \cdot |\ell|^{\frac{1}{4}} \left\| u^1_\ell \right\|_{L^2_t L^\infty_y}     \left\|  \theta_{k-\ell} \right\|_{L^\infty_t L^2_y} \mathrm{d} \ell \mathrm{d} k\\
			\lesssim& \nu^{-\frac{1}{2}} \int_{  |k-\ell| \ge 1} \int_{    |\ell| \ge 1}   |k|^{\frac{1}{12}}  |\ell|^{\frac{1}{4}} \cdot \nu^{-\frac{1}{8}} |\ell|^{-\frac{3}{4}} \cdot  |k-\ell|^{-\frac{1}{3}}   E_\ell[\omega_\ell]      \cdot  E_{k-\ell}[\theta_{k-\ell}]  \mathrm{d} \ell \mathrm{d} k\\
			\lesssim& \nu^{-\frac{1}{2}} \int_{  |k-\ell| \ge 1} \int_{    |\ell| \ge 1}     E_\ell[\omega_\ell]      \cdot  E_{k-\ell}[\theta_{k-\ell}]    \mathrm{d} \ell \mathrm{d} k \\
\lesssim& \nu^{-\frac{1}{2}} \mathcal{E}_{\omega} \cdot \mathcal{E}_{\theta}.
			\end{align*}
			In the penultimate step, we use the following observation:
			$$  |k|^{\frac{1}{12}}  |\ell|^{\frac{1}{4}} \cdot \nu^{-\frac{1}{8}} |\ell|^{-\frac{3}{4}} \cdot  |k-\ell|^{-\frac{1}{3}}  \lesssim \nu^{-\frac{1}{8}} |\ell|^{-\frac{1}{3}} |k-\ell|^{-\frac{1}{3}} |k|^{\frac{1}{12}} \lesssim \nu^{-\frac{1}{8}} |k|^{-\frac{1}{4}} \lesssim 1.$$
			By combining the estimates for all components in $R_{1,3}^{H,H}$, we obtain
			\begin{align*}
				R_{1,3}^{H,H} \lesssim \nu^{-\frac{1}{2}} \mathcal{E}_{\omega} \cdot \mathcal{E}_{\theta}.
			\end{align*}
			Furthermore, using the combined estimates for $R_{1,1}^{H,H}$ through $R_{1,3}^{H,H}$, we complete the proof of the lemma in the case
			$
			\nu^{\frac{1}{6}} |k|^{\frac{1}{3}} \ge 1.
			$
		\end{itemize}
	Finally, by synthesizing the discussions in {\bf Case 1} and {\bf Case 2}, we conclude the complete proof of the lemma.
	\end{proof}
\end{lemma}

\subsection{Completion the proof of the nonlinear stability of Theorem \ref{thm1.2}}
Next, we complete the proof of Theorem \ref{thm1.2}.
First, based on Lemmas \ref{prop_theta_low}--\ref{prop_theta_high}, which provide control over the nonlinear temperature terms, we can summarize
\begin{align}\label{ta_es}
	\mathcal{E}_{\theta} \lesssim \mathcal{E}^{\mathrm{in}}_{\theta}  + \nu^{-\frac{1}{2}} \mathcal{E}_{\omega} \cdot \mathcal{E}_{\theta}.
\end{align}
To this end, we define an energy functional as $$ \mathcal{E}_{\mathrm{total}} = \mathcal{E}_{\omega} + 2C \nu^{-\frac{1}{3}} \mathcal{E}_{\theta}.$$
Then, combining \eqref{ta_es} and \eqref{omega_es}, we obtain
$$ \mathcal{E}_{\mathrm{total}} \le  C \mathcal{E}^{\mathrm{in}}_{\mathrm{total}} + C \nu^{-\frac{1}{2}}  \mathcal{E}_{\mathrm{total}} ^2.$$
Using the initial conditions and choosing $\varepsilon_0$ and $\varepsilon_1$ sufficiently small, under the bootstrap hypothesis, we can deduce
\begin{align*}
	\mathcal{E}_{\mathrm{total}} \le C \nu^{\frac{1}{2}}.
\end{align*}
This completes the proof of Theorem \ref{thm1.2}. \hspace{7.8cm}
$\square$

\bigskip

\section*{Appendix}

\setcounter{theorem}{0}
\setcounter{equation}{0}
\renewcommand{\theremark}{A.\arabic{remark}}
\renewcommand{\theequation}{A.\arabic{equation}}
\renewcommand{\thetheorem}{A.\arabic{theorem}}
Let $H$ be a Hilbert space, and let $\mathcal{A}$ be a densely defined operator on $H$ with domain $D(\mathcal{A}) \subset H$.
We define the quantity
\[
\Phi(\mathcal{A}) := \inf \big\{ \| (\mathcal{A} - i \lambda) f \| : f \in D(\mathcal{A}), \, \lambda \in \mathbb{R}, \, \|f\| = 1 \big\},
\]
where $\|\cdot\|$ denotes the norm in $H$ and $i = \sqrt{-1}$. Intuitively, $\Phi(\mathcal{A})$ measures the distance of $\mathcal{A}$ from being purely imaginary in a certain sense.

\medskip

Now, we have the following lemma, which provides an upper bound for the semigroup generated by an $m$-accretive operator.
\begin{lemma}\label{le-GP}\cite{Wei2021SCM}
	Let $ \mathcal{A}$ be a $m$-accretive operator in a Hilbert space $ H$. Then $ \| e^{-\mathcal{A}t}\| \le e^{\frac{\pi}{2}}  e^{-\Phi(\mathcal{A})}$ for any $ t \ge 0$.
\end{lemma}
In other words, the decay of the semigroup is controlled by $\Phi(\mathcal{A})$, which depends on the resolvent behavior of $\mathcal{A}$ along the imaginary axis.

\begin{lemma}\label{lebasis1}
	Let $k,\ell \in \mathbb{R}$, and assume that $|k| \le 10\nu$ and $10\nu \le |\ell| \le 1$. If $10\nu \le |k-\ell| \le 2$, then
	\[
	|k-\ell| \approx |\ell|.
	\]
	\begin{proof}
		First, by the triangle inequality, we have
		\[
		|k-\ell| \le |k| + |\ell| \le 2|\ell|.
		\]
		
		On the other hand, we consider the lower bound:
		
		\begin{itemize}
			\item If $k,\ell > 0$, then $|k-\ell| = \ell - k \ge 10\nu$. Since $\ell \ge 10\nu + k \ge 2 k$, it follows that
			\[
			|k-\ell| = \ell - k \ge \frac{1}{2} |\ell|.
			\]
			
			\item If $k>0$ and $\ell < 0$, then $|k-\ell| = k - \ell \ge |\ell|$.
				
			The cases $k, \ell < 0$ and $k<0,\, \ell>0$ are completely analogous to the two cases above, and we omit the details.
		\end{itemize}
	
		Combining the upper and lower bounds, we conclude $|k-\ell| \approx |\ell|$.
		
	\end{proof}
\end{lemma}

\begin{lemma}\label{lebasis2}
	Let $k,\ell \in \mathbb{R}$ with $|k| \ge 1$ and $|\ell| \ge 1$. If $|k-\ell| \ge 1$, then for any $\alpha \ge 0$ and $m,n \in [0,\alpha]$, we have
	\[
	|k|^\alpha \lesssim |k|^m \cdot |k-\ell|^{\alpha-m} + |k|^n \cdot |\ell|^{\alpha-n}.
	\]
	\begin{proof}
		We consider two cases:
		
		\begin{itemize}
			\item If $|k-\ell| \le \frac{|k|}{2}$, then $|k| \approx |\ell|$, and the inequality holds trivially.
			
			\item If $|k-\ell| \ge \frac{|k|}{2}$, then the first term $|k|^m |k-\ell|^{\alpha-m} \gtrsim |k|^\alpha$, so the inequality is obvious.
		\end{itemize}
		
		Thus, the desired inequality is established in both cases.
	\end{proof}
\end{lemma}

\begin{lemma}\label{lebasis3}
	Let $\theta_k$ be a solution of system \eqref{equation_k}.
	Assume that $\nu^{\frac{1}{6}} |k|^{\frac{1}{3}} \ge 1$. Then
	\begin{align*}
	\| \theta_k \|_{L^\infty_t L^2_y}
	+ (\nu k^2)^{\frac{1}{6}} \| \theta_k \|_{L^2_t L^2_y}
	&\le
	\| \theta_k \|_{L^\infty_t L^2_y}
	+ (\nu k^2)^{\frac{1}{2}} \| \theta_k \|_{L^2_t L^2_y} \\
	&\le
	C \| \theta_k^{\mathrm{in}} \|_{L^2_y}
	+ C \nu^{-\frac{1}{2}}
	\big(
	\| g^1_k \|_{L^2_t L^2_y}
	+
	\| g^2_k \|_{L^2_t L^2_y}
	\big).
	\end{align*}
	
	\begin{proof}
		The first inequality follows immediately from the assumption
		\(
		\nu^{\frac{1}{6}} |k|^{\frac{1}{3}} \ge 1
		\),
		which implies
		\(
		(\nu k^2)^{\frac16} \le (\nu k^2)^{\frac12}.
		\)
		
		We now prove the second inequality.
		Taking the real part of the inner product of
		\(
		(\partial_t - \mathcal{A}) \theta_k
		\)
		with $\theta_k$ gives
		\[
		\operatorname{Re}
		\left\langle
		(\partial_t - \mathcal{A})\theta_k ,\, \theta_k
		\right\rangle
		=
		- \operatorname{Re}
		\left\langle
		ik g^1_k + \partial_y g^2_k,\, \theta_k
		\right\rangle .
		\]
		
		Using Young's inequality and integrating by parts in $y$, we obtain
		\begin{align*}
		\| \theta_k \|_{L^\infty_t L^2_y}^2
		+ \nu |k|^2 \| \theta_k \|_{L^2_t L^2_y}^2
		+ \nu \| \partial_y \theta_k \|_{L^2_t L^2_y}^2
		&\le
		\| \theta_k^{\mathrm{in}} \|_{L^2_y}^2
		+ \frac{1}{4} \nu
		\big(
		|k|^2 \| \theta_k \|_{L^2_t L^2_y}^2
		+
		\| \partial_y \theta_k \|_{L^2_t L^2_y}^2
		\big) \\
		&\quad
		+
		\nu^{-\frac{1}{2}}
		\big(
		\| g^2_k \|_{L^2_t L^2_y}^2
		+
		\| g^1_k \|_{L^2_t L^2_y}^2
		\big).
		\end{align*}
		
		Absorbing the $\frac14 \nu$ terms on the right-hand side completes the proof.
	\end{proof}
\end{lemma}

\vskip .2in
\section*{Acknowledgement}
\noindent{J. Wu was partially supported by the National Science Foundation of the United States under Grants DMS 2104682 and DMS 2309748. X. Zhai was partially supported by  the Guangdong Provincial Natural Science Foundation under grant 2024A1515030115. }

 \vskip .2in
\noindent{\bf Data Availability Statement} Data sharing is not applicable to this article as no
data sets were generated or analysed during the current study.

\vskip .2in

\noindent{\bf Conflict of Interest} The authors declare that they have no conflict of interest. The
authors also declare that this manuscript has not been previously published, and
will not be submitted elsewhere before your decision.


\begin{thebibliography}{10}
	
	\bibitem{Dolce2023CPAM} J. Bedrossian, R. Bianchini, M. Coti Zelati, M. Dolce, {Nonlinear inviscid damping and shear-buoyancy instability in the two-dimensional Boussinesq equations},	{\it Comm. Pure Appl. Math.}, {\bf 76}  (2023),  no. 12, 3685--3768.
	
\bibitem{BGM-BAMS} J. Bedrossian, P. Germain, N. Masmoudi,  Stability of the Couette flow at high Reynolds number in 2D and 3D, {\it Bull. Amer. Math. Soc.},  {\bf56}, 2019,  373--414.

\bibitem{Bm}
J. Bedrossian, N. Masmoudi,
Inviscid damping and the asymptotic stability of planar shear flows in the 2D Euler equations,
 {\it Publ. Math. Inst. Hautes \'Etudes Sci.}, {\bf122} (2015), 195--300.




	
	\bibitem{Dolce2021India} R. Bianchini, M. Coti Zelati, M. Dolce,{ Linear inviscid damping for shear flows near Couette in the 2D stably stratified regime}, {\it Indiana Univ. Math. J.}, {\bf 71}   (2022) 1467--1504.
	
	\bibitem{Cao2012} C. Cao, J. Wu, {Global regularity for the two-dimensional anisotropic Boussinesq equations with vertical dissipation}, {\it Arch. Ration. Mech. Anal.}, {\bf 208}  (2013) 985--1004.
	
	\bibitem{Chae2006} D. Chae, {Global regularity for the 2D Boussinesq equations with partial viscosity terms}, {\it Adv. Math.}, {\bf 203} (2006) 497--513.
	
	\bibitem{Chenqi2020ARMA} Q. Chen, T. Li, D. Wei,  Z. Zhang, {Transition threshold for the 2-D Couette flow in a finite channel}, {\it Arch. Ration. Mech. Anal.}, {\bf 238} (2020)  125--183.
	
	\bibitem{Chenqi2024AMS} Q. Chen, D. Wei,  Z. Zhang, {Transition threshold for the 3D Couette flow in a finite channel},  {\it Mem. Amer. Math. Soc.}, {\bf296} (2024), v+178 pp.
	
 \bibitem{Miao2025Arxiv} Q. Chen, Z. Li, C. Miao, {Quantitative stability for the 2D Couette flow on the infinite channel with non-slip boundary condition}, {\it arXiv:2510.18376}.

\bibitem{chenwangyang} Y. Chen, W. Wang, G. Yang, {Stability threshold of Couette flow for  Boussinesq equations in $\mathbb{R}^2$}, {\it arXiv:2508.11908}.




    \bibitem{ZZ2025}
		M. Coti Zelati and A. Del Zotto,
		Suppression of lift-up effect in the 3D Boussinesq equations around a stably stratified Couette flow,
		{\it Quart. Appl. Math.}, \textbf{83} (2025),  389--401.
		
		\bibitem{ZZW2024}
		M. Coti Zelati, A. Del Zotto and K. Widmayer,
		Stability of viscous three-dimensional stratified Couette flow via dispersion and mixing,  {\it arXiv:2402.15312V1}.




	\bibitem{Constantin} P. Constantin, C.R. Doering, {Heat transfer in convective turbulence}, {\it Nonlinearity}, {\bf 9}  (1996),  no. 4, 1049--1060.
	
\bibitem{CWW202501}
		S. Cui, L. Wang and W. Wang,
		Stability threshold of Couette flow for 3D Boussinesq system in Sobolev spaces,
		{\it  arXiv:2504.16401}.

	
	
	\bibitem{Danchin2011} R. Danchin, M. Paicu, {Global existence results for the anisotropic Boussinesq system in dimension two}, {\it Math. Models Methods Appl. Sci.}, {\bf 21}   (2011) 421--457.


{\bibitem{DHB} F. Daviaud, J. Hagseth and P. Berg\'e,  Subcritical transition to turbulence in plane Couette flow, {\it Phys. Rev. Lett.}, {\bf69}, 1992, 2511--2514.}


	
	\bibitem{Deng2021JFA} W. Deng, J. Wu, P. Zhang, {Stability of Couette flow for 2D Boussinesq system with vertical dissipation}, {\it J. Funct. Anal.},  {\bf 281}  (2021) 109255.
	
\bibitem{Dolce2024} M. Dolce. {Stability threshold of the 2D couette flow in a homogeneous magnetic field using symmetric variables}, {\it Comm. Math. Phys.}, {\bf 405}  (2024).

\bibitem{DR} P. Drazin, W. Reid, Hydrodynamic Stability, Cambridge Monographs Mech. Appl. Math., Cambridge
 Univ. Press, New York, 1981.


	\bibitem{Getling} A.V. Getling, {Rayleigh-B\'enard Convection, Advanced Series in Nonlinear Dynamics}, vol. 11, World Scientific Publishing Co., Inc., River Edge, NJ, 1998. Structures and dynamics.	
	
	\bibitem{Gill} A. Gill, {Atmosphere-Ocean Dynamics}, {\it International Geophysics}, vol. 30, Academic Press, Cambridge, 1982.
	
	\bibitem{Houyizhao} T.Y. Hou, C. Li, {Global well-posedness of the viscous Boussinesq equations}, {\it Discrete Contin. Dyn. Syst.} {\bf 12}   (2005) 1--12.

	\bibitem{Howard} L. N. Howard, {Note on a paper of john w. miles}, {\it J. Fluid Mech.} {\bf 10} (1961), 509--512.

\bibitem{Kel} L. Kelvin, Stability of fluid motion-rectilinear motion of viscous fluid between two
 parallel plates,  {\it  Phil. Mag.},  {\bf24}, 1887, 188--196.

\bibitem{Knobel} N. Knobel, {Suppression of fluid echoes and Sobolev stability threshold for 2D dissipative fluid equations around Couette flow}, {\it  arXiv:2505.23391}, 2025.


\bibitem{Kraichnan} R. H. Kraichnan, {Inertial Ranges in Two-dimensional Turbulence},{ \it  Phys. Fluids},  {\bf 10}  (1967)  1417--1423.

	\bibitem{Lai2011} M. Lai, R. Pan, K. Zhao, {Initial boundary value problem for two-dimensional viscous Boussinesq equations}, {\it Arch. Ration. Mech. Anal.}, {\bf 199}   (2011) 739--760.
	
	\bibitem{Lijinkai2016} J. Li, E.S. Titi, {Global well-posedness of the 2D Boussinesq equations with vertical dissipation}, {\it Arch. Ration. Mech. Anal.}, {\bf 220}   (2016) 983--1001.
	
\bibitem{liangtaoscm} T. Liang, J. Wu, X. Zhai, {Stability threshold analysis of the Boussinesq system with general viscosity near Couette flow}, {\it Preprint}.
	
\bibitem{liangtaofenshujie} T. Liang, J. Wu, X. Zhai, {Stability threshold of Couette flow for 2D Boussinesq equations  with fractional dissipation}, {\it Preprint}.

\bibitem{liangtaojimj} T. Liang, Y. Li, X. Zhai, {Stability threshold for two-dimensional Boussinesq system near-Couette shear flow in a finite channel}, {\it arXiv:2504.02669}.



	\bibitem{Majda}  A. Majda, {Introduction to PDEs and Waves for the Atmosphere and Ocean, Courant Lecture Notes in Mathematics}, vol. 9, New York University, Courant Institute of Mathematical Sciences/American Mathematical Society, New York/Providence, RI, 2003.
	
\bibitem{Miles} J. W. Miles, {On the stability of heterogeneous shear flows}, {\it J. Fluid Mech.} {\bf 10} (1961), 496--508.

	\bibitem{Majada2} A.J. Majda, A.L. Bertozzi,{ Vorticity and Incompressible Flow}, Cambridge Texts in Applied Mathematics, vol. 27, Cambridge University Press, Cambridge, 2002.
	
	\bibitem{Zhao2022ARMA} N. Masmoudi, B. Said-Houari, W. Zhao, {Stability of Couette flow for 2D Boussinesq system without thermal diffusivity}, {\it Arch. Ration. Mech. Anal.}, {\bf 245} (2022) 645--752.
	
	\bibitem{Masmoudi2023JFA} N. Masmoudi, C. Zhai,  W. Zhao, {Asymptotic stability for two-dimensional Boussinesq systems around the Couette flow in a finite channel}, {\it J. Funct. Anal.}, {\bf 284} (2023) 109736.
	
	
	
	\bibitem{Niu2024ARxiv} B. Niu, W. Zhao, { Improved stability threshold of the two-dimensional Couette flow for Navier-Stokes-Boussinesq systems via quasi-linearization}, {\it  arXiv:2409.16216}, 2024.
	




\bibitem{renxiaoxia2026} X. Ren, D. Wei,  On the stability threshold of Couette flow for 2D Boussinesq equations, {\it Nonlinear Anal. Real World Appl.}, {\bf87}, 2026, Paper No. 104421, 11 pp.


\bibitem{renxiaoxia}     X. Ren, D. Wei,  Transition threshold of Couette flow for 2D Boussinesq equations, arXiv:2506.03679 .


\bibitem{Re} O. Reynolds, An experimental investigation of the circumstances which determine whether the motion of water shall be direct or sinuous, and of the law of resistance in parallel channels, {\it Philos. Trans. R. Soc. Lond.}, {\bf174}, 1883, 935--982.



\bibitem{Rom}     V. A. Romanov,  Stability of plane-parallel Couette flow, {\it Funkcional. Anal. i Prilozen}, {\bf 7}, 1973, 62--73.






	\bibitem{Tao2017JDE} L. Tao, J. Wu, {The 2D Boussinesq equations with vertical dissipation and linear stability of shear flows}, {\it J. Differ. Equ.}, {\bf 267}  (2019) 1731--1747.
	
\bibitem{TA} N. Tillmark, P. H. Alfredsson,  Experiments on transition in plane Couette flow, {\it J. Fluid Mech.},
 {\bf235}, 1992, 89--102.

\bibitem{T} L. Trefethen, A. Trefethen, S. Reddy, T. Driscoll,  Hydrodynamic stability
 without eigenvalues, {\it Science}, {\bf261}, 1993, 578--584.



\bibitem{Wang2024arxiv}  F. Wang,  Z. Zhang. {The stability threshold for 2d mhd equations around couette with
	general viscosity and magnetic resistivity}, {\it arXiv:2410.20404}, 2024.
	
	\bibitem{wangwang} G. Wang, L. Wang, {On asymptotic stability of Couette flow for 2-D Boussinesq system in whole space via Green's function,}  {\it Discrete Contin. Dyn. Syst. Ser. B}, \textbf{30} (2025), no. 8, 3014--3041.

\bibitem{wangyang}
		W. Wang,  G. Yang,
		Stability threshold of the Couette flow for the 2D Boussinesq equations in a finite channel,
		preprint.
\bibitem{Wei2021SCM}
	D. Wei, Diffusion and mixing in fluid flow via the resolvent estimate, {\it Sci. China Math.}, \textbf{64} (2021), 507--518.

\bibitem{wei2}
		D. Wei,  Z. Zhang,
		{Transition threshold for the 3D Couette flow in Sobolev space},
		{\it Comm. Pure Appl. Math.}, \textbf{74} (2021), no. 11, 2398--2479.
		
	\bibitem{Wei2023Tunisi} D. Wei, Z. Zhang,  {Nonlinear enhanced dissipation and inviscid damping for the 2D Couette flow},  {\it Tunis. J. Math.}, {\bf 5}  ( 2023)  573--592.	
		
		\bibitem{wei4}
		D. Wei, Z. Zhang,
		{Stability threshold of the 2D Couette flow in a finite channel},
		preprint.
		
		\bibitem{wei5}
		D. Wei, Z. Zhang, W. Zhao,
		{Linear inviscid damping for a class of monotone shear flow in Sobolev spaces},
		 {\it Comm. Pure Appl. Math.}, \textbf{71} (2018), no. 4, 617--687.
		
		\bibitem{wei6}
		D. Wei, Z. Zhang, W. Zhao,
		{Linear inviscid damping and vorticity depletion for shear flows},
		 {\it Ann. PDE}, \textbf{5} (2019), no. 1, Paper No. 3, 101.
		
		\bibitem{wei7}
		D. Wei, Z. Zhang, W. Zhao,
		{Linear inviscid damping and enhanced dissipation for the Kolmogorov flow},
		 {\it Adv. Math.}, \textbf{362} (2020), 106963, 103 pp.
		
		\bibitem{wei8}
		D. Wei, Z. Zhang, H. Zhu,
		{Linear inviscid damping for the $\beta$-plane equation},
		 {\it Comm. Math. Phys.}, \textbf{375} (2020), no. 1, 127--174.





	\bibitem{Lin2018JMFM} J. Yang, Z. Lin, {Linear inviscid damping for Couette flow in stratified fluid}, {\it J. Math. Fluid Mech.}, {\bf 20}  (2018) 445--472.
	
	\bibitem{ZelatiARxiv} M. C. Zelati, A. D. Zotto, K. Widmayer, {Stability of viscous three-dimensional stratified
	Couette flow via dispersion and mixing}, {\it arXiv:2402.15312} (2024).

	\bibitem{Zhao2023SIAM} C. Zhai, W. Zhao, {Stability threshold of the Couette flow for Navier-Stokes Boussinesq system with large Richardson number $ \gamma^2 > \frac{1}{4}$ }, {\it SIAM J. Math. Anal.}, {\bf 55}   (2023) 1284--1318.
	
	\bibitem{Zhang2023JMPA} Z. Zhang, R. Zi, {Stability threshold of Couette flow for 2D Boussinesq
equations in Sobolev spaces}, {\it J. Math. Pures Appl.}, {\bf 179} (2023) 123--182.
	
\bibitem{Zillinger1}
C. Zillinger, On enhanced dissipation for the Boussinesq equations, {\it J. Differ. Equ.}, {\bf282} (2021) 407--445.

\bibitem{Zillinger2}
C. Zillinger, On the Boussinesq equations with non-monotone temperature profiles, {\it J. Nonlinear Sci.},  {\bf31} (2021) 64.

	


\end{thebibliography}
\end{document}